\documentclass[a4paper,10pt]{article}

   \usepackage{amsmath}
   \usepackage{amssymb}
   \usepackage{mathrsfs}
   \usepackage[dvips]{graphicx}
   \usepackage{amsthm}
   \usepackage[nice]{nicefrac}

   \newtheorem{theo}{Theorem}[section]
   \newtheorem{lemm}[theo]{Lemma}
   \newtheorem{prop}[theo]{Proposition}
   \newtheorem{rema}[theo]{Remark}
   \newtheorem{coro}[theo]{Corollary}
   \newtheorem{defi}[theo]{Definition}

   \def \calA{\mathcal A}
   \def \ha{\hat \alpha}

   \def \dtn{\dt^{(n)}}
   \def \dt{t_\delta}
   \def \on{\{0,\ldots,N\}}
   \def \hb{\hat \beta}
   
   \def \size{\hbox{\rm size}\,}
   \def \sgn{\text{sgn\,}}
   \def \onp{\{0,\ldots,N+1\}}
   \def \Pp{Problem $\mathcal P^k$}

   \def \n{\pmb{n}}
   \def \prt{\partial}
   
   \def \eps{\varepsilon}
   
   \def \N{\mathbb N}

   \newcommand{\dx}{\,\text dx}
   \newcommand{\dy}[1]{\,\text d#1}

\begin{document}

\begin{center}
   {\Large \bf Numerical approximation of a reaction-diffusion system \\
   \medskip with fast reversible reaction}
   \vskip0.2in

      {Robert Eymard \footnote{Universit\'e Paris-Est, 77454 Marne-la-Vall\'ee C\'edex 2, France; E-mail: {\tt Robert.Eymard@univ-mlv.fr}},
      \quad Danielle\ Hilhorst \footnote{Laboratoire de Math\'ematiques, CNRS and Universit\'e de Paris-Sud 11, 91405 Orsay C\'edex, France; E-mail: {\tt Danielle.Hilhorst@math.u-psud.fr}},
      \quad Michal Olech \footnote{Instytut Matematyczny Uniwersytetu Wroclawskiego, pl.\ Grunwaldzki 2/4, 50-384 Wroclaw, Polska; Laboratoire de Math\'ematiques, CNRS Universit\'e de Paris-Sud, 91405 Orsay C\'edex, France; E-mail: {\tt olech@math.uni.wroc.pl}. \\
      The preparation of this article has been partially supported by a Marie Curie Transfer of Knowledge Fellowship of the European Community's Sixth Framework Programme under contract number \texttt{MTKD-CT-2004-013389}}}
\end{center}

   \abstract{We consider the finite volume approximation of a reaction-diffusion system with fast reversible reaction. We deduce from {\it a priori} estimates that the approximate solution converges to the weak solution of the reaction\,--\,diffusion problem and satisfies estimates which do not depend on the chemical kinetics factor. It follows that the solution converges to the solution of a nonlinear diffusion problem, as the size of the volume elements and the time steps converge to zero while the kinetic rate tends to infinity.}

\noindent \textbf{Key words:}\ instantaneous reaction limit, mass-action kinetics, finite volume methods, convergence of approximate solutions, discrete \textit{a priori}\ estimates, Kolmogorov's theorem.

\noindent {\bf AMS subject classification:}\ 35K45, 35K50, 35K55, 65M12, 65N12, 65N22, 80A30, 92E20.

\section{Introduction}

In this paper, we consider chemical systems with fast reactions where mean reaction times vary from approximately $10^{-14}$ second to 1 minute. In particular, reactions that involve bond making or breaking are not likely to occur in less than $10^{-13}$ second. Moreover, chemical systems almost always involve some elementary reaction steps that are reversible and fast. \\
The study of reactions with rates that are outside of the time frame of ordinary laboratory operations requires specialized instrumentation, techniques and ways of proceeding (see for example Espenson \cite[Chapter 11]{espenson}). This work tries to give an efficient, quick and cheap way for numerical investigations of such reactions.

In this article, we consider a reversible chemical reaction between mobile species $\mathcal A$ and $\mathcal B$, that takes place inside a bounded region $\Omega \subset \mathbb R^d$ where $d = 1,2$ or $3$. If the region is isolated and diffusion is modelled by Fick's law, this leads to the reaction-diffusion system of partial differential equations
\begin{equation}
\begin{split}\label{FVMmain_eq}
   u_t &= a \Delta u - \alpha k \big(r_A(u) - r_B(v)\big) \text{\quad in \quad } \Omega \times (0,T), \\
   v_t &= b \Delta v + \beta  k \big(r_A(u) - r_B(v)\big) \text{\quad in \quad } \Omega \times (0,T),
\end{split}
\end{equation}
where $T > 0$ and $\Omega$ is a bounded set of $\mathbb R^d$. An example of explicit expressions and values for $\alpha,\ \beta,\ k,\ r_A,\ r_B,\ a,\ b$ is given in Section \ref{FVMnum_resul}. We supplement the system \eqref{FVMmain_eq} by the homogeneous Neumann boundary conditions
\begin{equation}\label{FVMmain_bd}
   \nabla u \cdot \n =\nabla v\cdot\n = 0 \text{\quad on \quad  } \prt \Omega \times (0,T),
\end{equation}
and the initial conditions of the form
\begin{equation}\label{FVMmain_in}
   u(x,0) = u_0(x),\ v(x,0) = v_0(x) \text{\quad in \quad } \Omega.
\end{equation}
In the sequel we call the system \eqref{FVMmain_eq} together with the boundary conditions \eqref{FVMmain_bd} and the initial conditions \eqref{FVMmain_in}, \Pp.\\
For a reversible reaction $m \mathcal A \rightleftharpoons n \mathcal B$ one has $\alpha = -m,\ \beta = n$ and the rate functions are of the form $r_A(u) = u^m$ and $r_B(v) = v^n$. Further discussion about this motivation and some concrete examples can be found in {\'E}rdi and T{\'o}th \cite{erdi_toth} and Espenson \cite{espenson}.

In practice, especially for ionic or radical reactions, changes due to reaction are often very fast compared to diffusive effects. This corresponds to a large rate constant $k$. Bothe and Hilhorst \cite{bothe_hilhorst} study the limit to an instantaneous reaction. They exploit a natural Lyapunov functional and use compactness arguments to prove that
\begin{equation*}
   u^k \rightarrow u \quad \text{and} \quad v^k \rightarrow v \quad \text{in} \quad L^2 \big(\Omega \times (0,T) \big),
\end{equation*}
as $k$ tends to infinity, where $(u^k,v^k)$ is the solution of \Pp\ and the limit $(u,v)$ is determined by
\begin{equation}\label{FVMpluis}
   r_A(u) = r_B(v) \quad \text{and} \quad \frac{u}{\alpha} + \frac{v}{\beta} = w,
\end{equation}
where $w$ is the unique weak solution of the nonlinear diffusion problem
\begin{equation}\label{FVMNonlDiff1}
\begin{split}
   w_t &= \Delta \phi(w) \quad \text{in} \quad \Omega \times (0,T)\\
   \frac{\prt \phi(w)}{\prt \n} &= 0 \quad \text{on} \quad \prt \Omega \times (0,T) \\
   w(x,0) & = w_0(x) := \frac{1}{\alpha}u_0(x) + \frac{1}{\beta}v_0(x) \quad \text{in} \quad \Omega \hskip 3pt ,
\end{split}
\end{equation}
with
\begin{equation}\label{FVMNonlDiff2}
\begin{split}
   \phi := \bigg(\frac{a}{\alpha} \text{id} + \frac{b}{\beta} \eta \bigg) &\circ \bigg( \frac{1}{\alpha} \text{id} + \frac{1}{\beta} \eta \bigg)^{-1} \quad \text{on} \quad \mathbb R^+, \\
   \eta &= r_B^{-1} \circ r_A.
\end{split}
\end{equation}
The identities in \eqref{FVMpluis} can be explained as follows:  the first one states that the system is in chemical equilibrium, while the second one defines $w$ as the quantity that is conserved under the chemical reaction. Given a function $w$, the system \eqref{FVMpluis} can be uniquely solved for $(u,v)$ if $r_A,\, r_B$ are strictly increasing with for instance $r_A(\mathbb R^+) \subset r_B(\mathbb R^+)$ so that $\eta = r_B^{-1} \circ r_A$ is well defined and strictly increasing. Under these assumptions $u$ is the unique solution of
\begin{equation*}
   \frac{1}{\alpha}u + \frac{1}{\beta}\eta(u) = w,
\end{equation*}
which gives the explicit representation of $u$ and $v$
\begin{equation}\label{FVMrepresUV}
   u = \bigg( \frac{1}{\alpha} \text{id} + \frac{1}{\beta} \eta \bigg)^{-1} \hskip -4pt \big( w \big),\quad v = \eta \circ \bigg( \frac{1}{\alpha} \text{id} + \frac{1}{\beta} \eta \bigg)^{-1} \hskip -4pt \big( w \big).
\end{equation}

\vskip 10pt
We assume the following hypotheses, which we denote by $\mathcal H$:
\textit{
\begin{enumerate}
   \item Let $\Omega$ be an open, connected and bounded subset of $\mathbb R^d$, where $d = 1,2$ or $3$, with a smooth boundary $\prt\Omega$,
   \item $u_0(x),\ v_0(x) \in L^\infty(\Omega)$ and there exist constants $U,\ V>0$ such that $0\leqslant u_0(x) \leqslant U$ and $0\leqslant v_0(x) \leqslant V$ in $\Omega$,
   \item $\alpha$, $\beta$, $a$, $b$ and $k$ are strictly positive real values (sometimes we use the notation $k \alpha = \ha$ and $k\beta = \hb$),
   \item Let $r_A(x),\ r_B(x)\in C^1(\mathbb R)$ be strictly increasing functions, such that $r_A(0) = r_B(0) = 0$, and assume furthermore that $r_A(\mathbb R^+) \subset r_B(\mathbb R^+)$.
\end{enumerate}
}

We recall from Bothe and Hilhorst \cite[Section 2]{bothe_hilhorst} that \Pp\ has a unique classical solution $(u^k,v^k)$ on every finite time interval $[0,T]$, for all nonnegative bounded initial data. By classical solution, we mean a function pair $(u^k,\,v^k)$ such that $u^k,\,v^k \in C^{2,1}\big( \Omega \times (0,T] \big) \cap C^{1,0}\big(\overline \Omega \times (0,T] \big)$ with $u^k, v^k \in C\big([0,T];\, L^2(\Omega) \big)$ (see also Lady{\v{z}}enskaja, Solonnikov and Ural{'c}eva \cite{lady_solo_ural}).

Next we present a notion of a weak solution of \Pp, which will be used in the sections \ref{FVMfixed_k} and  \ref{FVMunifcase}.
\begin{defi}\label{FVMweakdefi}
   We say that $(u^k,v^k)$ is a weak solution to \Pp\ if and only if
   \begin{enumerate}
      \item $u^k, v^k\in L^2\big(0,T;H^1(\Omega)\big)$ and $u^k_t, v^k_t\in L^2\big(0,T;(H^1(\Omega))'\big)$;
      \item Let $\Psi$ be the set of test functions, defined as
      \begin{equation*}
         \Psi = \Big \{ \psi \in C^{2,1}(\overline \Omega \times [0,T]): \nabla \psi \cdot \n = 0 \text{\quad on \quad} \partial \Omega \times [0,T] \text{\quad and \quad} \psi(T) = 0 \Big \}.
      \end{equation*}
      For a.e. $t\in(0,T)$ and all $\psi \in \Psi$
      \begin{multline}\label{FVMweakdef1}
         \int_\Omega u_0(x) \psi(x,0) \dx + \int_\Omega u^k(x,t) \psi_t(x,t) \dx + a \int_\Omega u^k(x,t) \Delta \psi(x,t) \dx \\
         - \hat \alpha \int_\Omega \psi(x,t) \Big(r_A \big(u^k(x,t) \big) - r_B \big(v^k(x,t)\big)\Big) \dx = 0
      \end{multline}
      and
      \begin{multline}\label{FVMweakdef2}
         \int_\Omega v_0(x) \psi(x,0) \dx + \int_\Omega v^k(x,t) \psi_t(x,t) \dx + b \int_\Omega u^k(x,t) \Delta \psi(x,t) \dx\\
         + \hat \beta \int_\Omega \psi(x,t) \Big(r_A \big(u^k(x,t) \big) - r_B \big(v^k(x,t) \big)\Big) \dx = 0.
      \end{multline}
   \end{enumerate}
\end{defi}
We remark that every essentially bounded weak solution of \Pp, in the sense of Definition \ref{FVMweakdefi}, is also a classical solution.

This paper is organized as follows. In section 2 we define a finite volume discretisation and an approximate solution $(u_{\cal D}^k,v_{\cal D}^k)$ for \Pp. In section 3 we prove a discrete comparison principle which yields  discrete $L^\infty$ estimates, and we show the existence and uniqueness of the approximate solution. Section 4 contains technical lemmas used further in the convergence proofs. The convergence of the approximate solution to the classical solution of \Pp\ in the case of fixed $k$ is proved in section 5. In section 6 we use a suitable Lyapunov function and we obtain a discrete $L^2\big(0,T;H^1(\Omega)\big)$ estimate, which does not depend on $k$. We then apply Kolmogorov's theorem and deduce the convergence of the approximate solutions to the classical solution of \Pp. Afterwards we show that the approximate solution $(u_{\cal D}^k,v_{\cal D}^k)$ converges to $(u,v)$ defined in \eqref{FVMrepresUV} as $k$ tends to $\infty$ and the size of the discretisation parameters tends to zero.\\
In Section \ref{FVMnum_resul} we present numerical results obtained with our finite volume scheme, for the reversible dimerisation of $o$-phenylene\-dio\-xydi\-me\-thyl\-si\-la\-ne ($2,2$-dimethyl-$1,2,3$-benzodioxasilole) which is a reaction of the type $2 \mathcal A \rightleftharpoons \mathcal B$ (see Meyer, Klein and Weiss \cite{meyer_klein_weiss}). On the one hand we compute the approximate solution $(u_{\cal D}^k,v_{\cal D}^k)$ of the solution $(u^k,v^k)$ of \Pp\ and on the other hand the numerical approximation $w_{\cal D}$ of the solution $w$ of the problem \eqref{FVMNonlDiff1}\,--\,\eqref{FVMNonlDiff2}, and we check that
\begin{equation*}
   \rule{20pt}{0pt} u_{\cal D}^k \approx \bigg( \frac{1}{\alpha} \text{id} + \frac{1}{\beta} \eta \bigg)^{-1} \hskip -4pt \big( w_{\cal D} \big)\quad \text{and} \quad v_{\cal D}^k \approx \eta \circ \bigg( \frac{1}{\alpha} \text{id} + \frac{1}{\beta} \eta \bigg)^{-1} \hskip -4pt \big( w_{\cal D} \big) \rule{20pt}{0pt}
\end{equation*}
for $k$ large enough and $\size({\cal D})$ small enough.

\begin{rema}
In what follows we denote by $C$, $C_k$ and $C_\psi$ positive generic constants which may vary from line to line.
\end{rema}

\section{The finite volume scheme}

The finite volume method has first been developed by engineers in order to study complex, coupled physical problems where the conservation of quantities such as masses, energy or impulsion must be carefully respected by the approximate solution. Another advantage of this method is that a large variety of meshes can be used in the computations. The finite volume methods are particularly well suited for numerical investigations of conservations laws. They are one of the most popular methods among the engineers performing computations for industrial purposes: the modelling of flows in porous media, problems related to oil recovery, questions related to hydrology, such as the numerical approximation of a stationary incompressible Navier\,--\,Stokes equations.

For a comprehensive discussion about the finite volume method, we refer to Eymard, Gallou{\"e}t and Herbin \cite{ciarlet_lions_VII} and the references therein.

\vskip 10pt
Following \cite{ciarlet_lions_VII}, we define a finite volume  discretization of $Q_T$.

\begin{defi}[Admissible mesh of $\Omega$]\label{FVMadmmesh}
   An admissible mesh ${\cal M}$ of $\Omega$ is given by a set of open, bounded subsets of\ $\Omega$ (control volumes) and a family of points (one per control volume), satisfying the following properties
   \begin{enumerate}
      \item The closure of the union of all the control volumes is $\overline \Omega$. We denote by $m_K$ the measure of each volume element $K$ and
      \begin{equation*}
         \size({\cal M})=\max_{K\in{\cal M}} m_k.
      \end{equation*}
      \item $K\cap L = \emptyset$ for any $(K,L) \in {\cal M}^2$, such that $K\neq L$. If $\overline{K} \cap \overline{L}\neq\emptyset$, then it is a subset of a hyperplane in $\mathbb R^d$. Let us denote by ${\cal E}\subset{\cal T}^2$ the set of pairs $(K,L)$, such that the $d-1$ Lebesgue measure of $\overline{K} \cap \overline{L}$ is strictly positive. For $(K,L)\in{\cal E}$ we write $K|L$ for the set $\overline{K} \cap \overline{L}$ and $m_{K|L}$ for the $d-1$ Lebesgue measure of $K|L$.
      \item For any $K \in {\cal M}$ we also define ${\cal N}_K = \{L\in{\cal T},\ (K,L)\in{\cal E}\}$ and assume that $\partial K=\overline{K} \backslash K= \big(\overline{K}\cap\prt\Omega\big) \cup \Big( \bigcup_{L \in {\cal N}_K} K|L \Big)$.
      \item There exists a family of points $(x_K)_{K \in {\cal M}}$, such that $x_K \in K$ and if $L\in{\cal N}_K$ then the straight line $(x_K,x_L)$ is orthogonal to $K|L$. We set
      \begin{equation*}
         d_{K|L}=d(x_K,x_L) \quad \text{and} \quad  T_{K|L}=\frac{m_{K|L}}{d_{K|L}},
      \end{equation*}
      where the last quantity is sometimes called the \textsl{transmissibility} across the edge $K|L$.
   \end{enumerate}
\end{defi}\medskip

\hskip -\parindent Since Problem \Pp is a time evolution problem, we also need to discretize the time interval $(0,T).$
\begin{defi}[Time discretization]\label{FVMadmtime}
   A time discretization of the interval $(0,T)$ is given by an integer value $N$ and by a strictly increasing sequence of real values $(t^{(n)})_{n\in\onp}$ with $t^{(0)} = 0$ and $t^{(N+1)} = T$. The time steps are defined by
   \begin{equation*}
      \dt^{(n)} = t^{(n+1)}-t^{(n)}\quad\text{for}\quad n\in\on.
   \end{equation*}
\end{defi}\medskip

\hskip -\parindent We may then define a discretization of the whole domain $Q_T$ in the following way.
\begin{defi}[Discretization of $Q_T$]\label{FVMdiscretization}
   A finite volume discretization ${\cal D}$ of $Q_T$ is defined as
   \begin{equation*}
      {\cal D} = \Big({\cal M}, {\cal E}, (x_K)_{K \in {\cal T}},  (t^{(n)})_{n\in\onp}\Big),
   \end{equation*}
where ${\cal M}$, ${\cal E}$ and $(x_K)_{K\in{\cal T}}$ are given in Definition
\ref{FVMadmmesh} and the sequence $(t^{(n)})_{n\in\onp}$ is a time discretization of $(0,T)$ in the sense of Definition \ref{FVMadmtime}. One then sets
   \begin{equation*}
      \size({\cal D}) = \max \Big \{ \size ({\cal M}) ,\dt^{(n)}:\ {n\in\on} \Big \}.
   \end{equation*}
\end{defi}\medskip

We present below the finite volume scheme which we use and define approximate solutions. We assume that the hypotheses $\mathcal H$ hold and suppose that ${\cal D}$ be an admissible discretization of $Q_T$ in the sense of Definition \ref{FVMdiscretization}. We prescribe the approximate initial conditions
\begin{equation}
   u_K^{(0)} = \frac 1 {m_K} \int_K u_0(x) \dx \quad \text{and}\quad
   v_K^{(0)} = \frac 1 {m_K} \int_K v_0(x) \dx,
\label{FVMschini}
\end{equation}
where $K\in{\cal M}$, and associate to \Pp\ the finite volume scheme
\begin{equation}
\begin{split}\label{FVMscheme}
   m_K \big(u_K^{(n+1)} &- u_K^{(n)}\big) - \dt^{(n)} a \sum_{L\in{\cal N}_K} T_{K|L} (u_L^{(n+1)} - u_K^{(n+1)}) + \dt^{(n)} m_K \ha \Big(r_A \big(u_K^{(n+1)} \big) - r_B \big( v_K^{(n+1)} \big)\Big) = 0,\\
   m_K \big(v_K^{(n+1)} &- v_K^{(n)}\big) - \dt^{(n)} b \sum_{L\in{\cal N}_K} T_{K|L} (v_L^{(n+1)} - v_K^{(n+1)}) - \dt^{(n)} m_K \hb \Big(r_A \big(u_K^{(n+1)} \big) - r_B \big( v_K^{(n+1)}\big) \Big) = 0.
\end{split}
\end{equation}
Note that \eqref{FVMscheme} is a nonlinear system of equations in the unknowns \begin{equation*}
   (u_K^{(n+1)},v_K^{(n+1)})_{K\in{\cal M},\ n\in\on}.
\end{equation*}
For $x\in\Omega$ and $t\in(0,T)$ let $K\in\cal M$ be such that $x\in K$ and $n\in\on$ be such that $t^{(0)} = 0$, $t^{(N+1)} = T$ and $t\in \big(t^{(n)},t^{(n+1)}\big]$. We can then define the approximate solutions
\begin{equation}\label{FVMuapp}
   u_{\cal D}(x,t) = u_K^{(n+1)} \mbox{\quad and \quad} v_{\cal D}(x,t) = v_K^{(n+1)} \hskip 5pt .
\end{equation}

\noindent In the next section, we will prove the existence and uniqueness of the solution of the discrete problem \eqref{FVMscheme}, together with the initial values \eqref{FVMschini}.

\section{The approximate solution}

In this section we prove the existence and uniqueness of the solution of the system \eqref{FVMscheme}. Let us start with a discrete version of the comparison principle.
\begin{prop}[Discrete comparison principle]\label{FVMcomprinciple}
   We suppose that the hypotheses $\mathcal H$ are satisfied. Let $\cal D$ be a discretization as in Definition \ref{FVMdiscretization}. Let  $(u_K^{(0)},v_K^{(0)})_{K\in{\cal M}}$ and $(\tilde u_K^{(0)},\tilde v_K^{(0)})_{K\in{\cal M}}$ be given sequences of real values such that
   \begin{equation*}
      u_K^{(0)} \leqslant \tilde u_K^{(0)}\quad \text{and} \quad v_K^{(0)} \leqslant \tilde v_K^{(0)},
   \end{equation*}
   for all $K\in{\cal M}$. If the sequences $(u_K^{(n+1)},v_K^{(n+1)})_{K\in{\cal M},\ n\in\on}$ and \\
   $(\tilde u_K^{(n+1)},\tilde v_K^{(n+1)})_{K\in{\cal M},\ n\in\on}$ satisfy the equations \eqref{FVMscheme} with the initial va\-lues $(u_K^{(0)},v_K^{(0)})_{K\in{\cal M}}$ and $(\tilde u_K^{(0)},\tilde v_K^{(0)})_{K\in{\cal M}}$, respectively, then for $K\in{\cal M}$ and $n\in\on$
   \begin{equation}\label{FVMeqcomprin}
      u_K^{(n+1)} \leqslant \tilde u_K^{(n+1)}\quad \text{and} \quad v_K^{(n+1)} \leqslant \tilde
      v_K^{(n+1)}.
   \end{equation}
\end{prop}
{\hskip -\parindent \bf Proof \hskip 5pt}
   We set $\hat u_K^{(n)} = u_K^{(n)} - \tilde u_K^{(n)}$ and $\hat v_K^{(n)} = v_K^{(n)} - \tilde v_K^{(n)}$ for all $K\in{\cal M}$ and $n\in\onp$ and define
   \begin{equation*}
   \begin{split}
      \hat A_K^{(n+1)} =& \Big( r_A \big(u_K^{(n+1)}\big) - r_A \big(\tilde u_K^{(n+1)}\big)\Big)/\hat u_K^{(n+1)},\\
      \hat B_K^{(n+1)} =& \Big( r_B \big(v_K^{(n+1)}\big) - r_B \big(\tilde v_K^{(n+1)}\big)\Big)/\hat v_K^{(n+1)},
   \end{split}
   \end{equation*}
   whenever $\hat u_K^{(n+1)}\neq 0$ (else $\hat A_K^{(n+1)} = 0$) or $\hat v_K^{(n+1)}\neq 0$ (else $\hat B_K^{(n+1)} = 0$). Since the functions $r_A$ and $r_B$ are monotone increasing, it follows that  $\hat A_K^{(n+1)}$ and $\hat B_K^{(n+1)}$ are nonnegative. We then have, by subtracting the discrete equation \eqref{FVMscheme} for $u_K^{(n+1)}$ and for ${\tilde u_K}^{(n+1)}$,
   \begin{equation}\label{FVMafter_ins}
   \begin{split}
      m_K \Bigg( 1 + \dt^{(n)} & \bigg(\hat \alpha \hat A_K^{(n+1)} + \frac{a}{m_K} \sum_{L \in{\cal N}_K} T_{K|L}\bigg) \Bigg) \hat u_K^{(n+1)}\\
      = &\ m_K \hat u_K^{(n)} + \dt^{(n)} a \sum_{L\in{\cal N}_K} T_{K|L} \hat u_L^{(n+1)}\\
      &\ + \dt^{(n)} m_K \hat\alpha \Big(r_B \big(v_K^{(n+1)} \big) - r_B \big(\tilde v_K^{(n+1)} \big  )\Big),
   \end{split}
   \end{equation}
   for  $K\in{\cal M}$ and $n\in\on$. Setting $s^+ = \max(s,0)$ and using that $s\leqslant s^+,\ (s+t)^+\leqslant s^+ + t^+$ we obtain
   \begin{equation}\label{FVMpart1}
   \begin{split}
      m_K  \Bigg( 1 + \dt^{(n)} \bigg(& \hat \alpha \hat A_K^{(n+1)} + \tfrac{a}{m_K} \sum_{L\in{\cal N}_K} T_{K|L}\bigg) \Bigg) \hat u_K^{(n+1)}\\
      \leqslant &\ m_K \big(\hat u_K^{(n)}\big)^+ +  \dt^{(n)} a \sum_{L\in{\cal N}_K} T_{K|L} \big(\hat u_L^{(n+1)}\big)^+ + \dt^{(n)} m_K \hat \alpha \Big(r_B \big(v_K^{(n+1)} \big) - r_B \big( \tilde v_K^{(n+1)} \big)\Big)^+,
   \end{split}
   \end{equation}
   where $K\in{\cal M}$ and $n\in\on$. Next we multiply the inequality \eqref{FVMpart1} by indicator of the set where $\hat u_K^{(n+1)}$ is nonnegative. Since the right-hand-side of \eqref{FVMpart1} is nonnegative as well, we obtain, acting similarly for both components,
   \begin{equation}\label{FVMcomp2}
   \begin{split}
      m_K \Bigg( 1 + \dt^{(n)}\bigg(& \hat \alpha \hat A_K^{(n+1)} + \frac{a}{m_K} \sum_{L\in{\cal N}_K} T_{K|L}\bigg) \Bigg) \big(\hat u_K^{(n+1)}\big)^+\\
      \leqslant &\ m_K \big(\hat u_K^{(n)}\big)^+ + \dt^{(n)} a \sum_{L \in{\cal N}_K} T_{K|L} \big(\hat u_L^{(n+1)}\big)^+ + \dt^{(n)} m_K \hat \alpha \Big(r_B \big(v_K^{(n+1)}\big) - r_B \big( \tilde v_K^{(n+1)}\big)\Big)^+,\\
      m_K \Bigg( 1 + \dt^{(n)}\bigg(& \hat \beta \hat B_K^{(n+1)} + \frac{b}{m_K} \sum_{L\in{\cal N}_K} T_{K|L}\bigg) \Bigg) \big(\hat v_K^{(n+1)} \big)^+\\
      \leqslant &\ m_K \big(\hat v_K^{(n)}\big)^+ + \dt^{(n)} b \sum_{L \in{\cal N}_K} T_{K|L} \big(\hat v_L^{(n+1)}\big)^+ + \dt^{(n)} m_K \hat \beta \Big(r_A \big(u_K^{(n+1)}\big) - r_A \big(\tilde u_K^{(n+1)} \big) \Big)^+.
   \end{split}
   \end{equation}
   Since
   \begin{equation*}
   \begin{split}
      \hat A_K^{(n+1)} (\hat u_K^{(n+1)})^+ = & \big(r_A(u_K^{(n+1)}) - r_A(\tilde u_K^{(n+1)})\big)^+ \\
      \hat B_K^{(n+1)} (\hat v_K^{(n+1)})^+ = & \big(r_B(v_K^{(n+1)}) - r_B(\tilde v_K^{(n+1)})\big)^+
   \end{split}
   \end{equation*}
   we add the first equation of \eqref{FVMcomp2} divided by $\hat\alpha$ and the
   second equation of \eqref{FVMcomp2} divided by $\hat\beta$, which yields
   \begin{multline}\label{FVMtech01}
      m_K \bigg( \frac{1}{\hat \alpha} + \dt^{(n)} \frac a {m_K \hat \alpha}  \sum_{L \in{\cal N}_K} T_{K|L}\bigg) \big(\hat u_K^{(n+1)}\big)^+ +\ m_K \bigg( \frac 1 {\hat \beta} + \dt^{(n)} \frac b {m_K \hat \beta} \sum_{L \in{\cal N}_K} T_{K|L}\bigg) \big(\hat v_K^{(n+1)}\big)^+\\
      \leqslant m_K \frac{1}{\hat \alpha} \big(\hat u_K^{(n)}\big)^+ + \dt^{(n)} \frac{a}{\hat \alpha} \sum_{L \in{\cal N}_K} T_{K|L} \big(\hat u_L^{(n+1)}\big)^+ +\ m_K \frac{1}{\hat \beta}\big(\hat v_K^{(n)}\big)^+ + \dt^{(n)} \frac{b}{\hat \beta} \sum_{L \in{\cal N}_K} T_{K|L} \big(\hat v_L^{(n+1)}\big)^+,
   \end{multline}
   for $K\in{\cal M}$ and $n\in\on$. Let us note that
   \begin{multline*}
      \rule{40pt}{0pt} \sum_{K \in {\cal M}}\sum_{L \in{\cal N}_K} T_{K|L} \big(\hat u_K^{(n+1)}\big)^+ = \sum_{L \in{\cal M}}\sum_{K \in{\cal N}_L} T_{L|K} \big( \hat u_L^{(n+1)} \big)^+\\
      = \sum_{L \in{\cal M}}\sum_{K \in{\cal N}_L} T_{K|L} \big(\hat u_L^{(n+1)} \big)^+ = \sum_{K \in{\cal M}}\sum_{L \in{\cal N}_K} T_{K|L} \big(\hat u_L^{(n+1)} \big)^+.\rule{40pt}{0pt}
   \end{multline*}
   Summing the inequalities \eqref{FVMtech01} over $K\in{\cal M}$, we get
   \begin{equation*}
      \sum_{K \in{\cal M}} \bigg[m_K \Big(\frac{1}{\hat \alpha} \big(\hat u_K^{(n+1)} \big)^+ + \frac{1}{\hat \beta} \big(\hat v_K^{(n+1)} \big)^+ \Big) \bigg] \leqslant \sum_{K \in{\cal M}} \bigg[ m_K \Big(\frac{1}{\hat \alpha} \big(\hat u_K^{(n)}\big)^+ + \frac{1}{\hat \beta} \big(\hat v_K^{(n)}\big)^+ \Big) \bigg],
   \end{equation*}
   which therefore leads, by induction, to
   \begin{equation*}
      \sum_{K \in{\cal M}} \bigg[m_K \Big(\frac{1}{\hat \alpha} \big(\hat u_K^{(n+1)} \big)^+ +\frac{1}{\hat\beta} \big(\hat v_K^{(n+1)} \big)^+ \Big) \bigg] = 0,
   \end{equation*}
   where $n\in\on$. It implies that $(\hat u_K^{(n+1)})^+=(\hat v_K^{(n+1)})^+=0$, which completes the proof. {\hskip 5pt $\blacksquare$ \medskip}

\begin{coro}[Discrete contraction in $L^1$ property]
   With the notation from Proposition \ref{FVMcomprinciple}, we have that
   \begin{multline*}
      \sum_{K\in \cal M} m_K  \bigg( \frac{\big|u_k^{(n+1)}-\tilde u_k^{(n+1)}\big|}{\hat\alpha} + \frac{\big|v_k^{(n+1)}-\tilde v_k^{(n+1)}\big|}{\hat\beta}\bigg) \leqslant \ \sum_{K\in \cal M} m_K  \bigg( \frac{\big|u_k^{(n)}-\tilde u_k^{(n)}\big|}{\hat\alpha} + \frac{\big|v_k^{(n)}-\tilde v_k^{(n)}\big|}{\hat\beta}\bigg)
   \end{multline*}
   for $n \in \on$. In other words, the discrete counterpart of the $L^1(\Omega)$-contraction property for solutions of \eqref{FVMmain_eq} (see e.g.\ \cite{bothe_hilhorst}) is preserved by the numerical scheme \eqref{FVMscheme}.
\end{coro}
{\hskip -\parindent \bf Proof \hskip 5pt}
   The proof directly follows from the proof of Proposition \ref{FVMcomprinciple}. Let us consider the term $\hat u_K$. We multiply the equation \eqref{FVMafter_ins} by $\sgn \big(\hat u_K^{(n+1)}\big)$. Then, the inequality $x \leqslant |x|$ yields
   \begin{equation*}
   \begin{split}
      m_K\Bigg( 1 + \dt^{(n)}\bigg(& \hat\alpha \hat A_K^{(n+1)} + a \sum_{L\in{\cal N}_K} T_{K|L}\bigg) \Bigg) \big|\hat u_K^{(n+1)}\big|\\
      \leqslant &\ m_K \big|\hat u_K^{(n)}\big| + \dt^{(n)} \tfrac{a}{m_K} \sum_{L\in{\cal N}_K} T_{K|L} \big|\hat u_L^{(n+1)}\big| + \dt^{(n)} m_K \hat\alpha \Big|r_B(v_K^{(n+1)}) - r_B(\tilde v_K^{(n+1)})\Big|.
   \end{split}
   \end{equation*}
   We proceed in the same way for $\hat v_K^{(n+1)}$ and remark that
   \begin{equation*}
   \begin{split}
      \hat A_K^{(n+1)} \big|\hat u_K^{(n+1)}\big| = & \Big|r_A(u_K^{(n+1)}) - r_A(\tilde u_K^{(n+1)})\Big|,\\
      \hat B_K^{(n+1)} \big|\hat v_K^{(n+1)}\big| = & \Big|r_B(v_K^{(n+1)}) - r_B(\tilde v_K^{(n+1)})\Big|,
   \end{split}
   \end{equation*}
   which enable us to obtain the counterpart of the inequalities in \eqref{FVMtech01} which we sum over $K\in\cal M$, as in the proof of Proposition \ref{FVMcomprinciple}. This yields the result. {\hskip 5pt $\blacksquare$ \medskip}

We now are in a position to prove a discrete $L^\infty$ estimate for the approximate solution.

\begin{theo}\label{FVMlinf}
   Let ${\cal D} = \Big({\cal M}, {\cal P}, {\cal E}, (\dt^{(n)})_{n \in\onp} \Big)$ be an admissible discretization of $Q_T$ in the sense of Definition \ref{FVMdiscretization}. We suppose that the hypotheses $\mathcal H$ are satisfied. Let  $(u_K^{(0)},v_K^{(0)})_{K\in{\cal M}}$ be given by \eqref{FVMschini} and $(u_K^{(n+1)},v_K^{(n+1)})$ satisfy \eqref{FVMscheme} for $K\in{\cal M}$ and $n\in\on$. Then
   \begin{equation}\label{FVMeqlinf}
      0 \leqslant u_K^{(n+1)} \leqslant U+ \frac \alpha \beta V \quad \text{and}\quad 0 \leqslant v_K^{(n+1)}
      \leqslant V+ \frac\beta  \alpha U,\
   \end{equation}
   for all $K\in{\cal M}$ and $\ n\in\on$, where $U$ and $V$ are the positive constants from the hypothesis $\mathcal H\, 2$.
\end{theo}
{\hskip -\parindent \bf Proof \hskip 5pt}
   From Proposition \ref{FVMcomprinciple} we immediately obtain that $u_K^{(n+1)}$ and $v_K^{(n+1)}$ are nonnegative for $K\in{\cal M}$ and $\ n\in\on$. In order to find a discrete upper solution, we consider approximate solutions of the corresponding system of ordinary differential equations. More precisely, we consider sequences $(\bar u^{(n)})_{n\in\onp}$, $(\bar v^{(n)})_{n\in\onp}$ (we postpone for a moment the proof that they exist) such that
   \begin{equation*}
      \bar u^{(0)} = U, \quad \bar v^{(0)} = V
   \end{equation*}
   and
   \begin{equation}\label{FVMunisch}
   \begin{split}
      \bar u^{(n+1)} - \bar u^{(n)} &= \alpha k\,\dt^{(n)} \Big(r_B(\bar v^{(n+1)}\big) - r_A \big(\bar
         u^{(n+1)} )\Big),\\
      \bar v^{(n+1)} - \bar v^{(n)} &= \beta k\, \dt^{(n)} \Big(r_A \big(\bar u^{(n+1)} ) - r_B(\bar v^{(n+1)} \big)\Big),
   \end{split}
   \end{equation}
   for $n\in\on$. We note that the sequences $(\bar u^{(n+1)})_{n\in\on}$, $(\bar v^{(n+1)})_{n\in\on}$ satisfy \eqref{FVMscheme} with the initial data $U,V$. Therefore they satisfy the comparison principle from Proposition \ref{FVMcomprinciple} which yields
   \begin{equation}\label{FVMunischpos}
      0 \leqslant \bar u^{(n+1)} \quad \text{and} \quad 0 \leqslant \bar v^{(n+1)}\quad \text{for all} \quad n \in \on.
   \end{equation}
   Adding up the first equation of \eqref{FVMunisch} divided by $\alpha$ and the second one divided by $\beta$, we obtain
   \begin{equation*}
      \frac{\bar u^{(n+1)}}{\alpha}  + \frac{\bar v^{(n+1)}}{\beta}
    = \frac{\bar u^{(n)}}{\alpha}  + \frac{\bar v^{(n)}}{\beta}
    = \ldots = \frac{U}{\alpha} + \frac{V}{\beta} \hskip 5pt.
   \end{equation*}
   We deduce from the previous equation and from \eqref{FVMunischpos} that
   \begin{equation}\label{FVMequnilinf}
      0 \leqslant \bar u^{(n+1)} \leqslant U+ \frac \alpha \beta V\text{\quad and\quad } 0 \leqslant \bar
      v^{(n+1)}
      \leqslant V+ \frac\beta  \alpha U,
   \end{equation}
   for $n\in\on$. {\hskip 5pt $\blacksquare$ \medskip}

In order to prove the existence of the sequences $(\bar u^{(n)})_{n\in\onp}$ and $(\bar v^{(n)})_{n\in\onp}$ we use the topological degree theory in finite dimensional spaces. The reader can find basic definitions as well as further informations about this powerful theory in Deimling \cite{deimling}. An example of the application of this tool to the analysis of finite volume schemes can be found in Eymard, Gallou{\"e}t, Ghilani and Herbin \cite{re_tg_ghilani_re}.

With $\mathcal F,\mathcal G: \mathbb R^2\rightarrow\mathbb R^2$ defined as
\begin{equation*}
\begin{split}
   \mathcal F \big ( \bar u^{(n)}, \bar v^{(n)} \big ) =&\ \big ( \bar u^{(n)}, \bar v^{(n)} \big ),\\
   \mathcal G \big ( \bar u^{(n)}, \bar v^{(n)} \big ) =&\ \bigg ( \hat \alpha \dt^{(n-1)} \Big(r_A \big(\bar u^{(n)} ) - r_B(\bar v^{(n)}\big)\Big),\\
   & -\hat\beta \dt^{(n-1)} \Big(r_A \big(\bar u^{(n)} ) - r_B(\bar v^{(n)} \big)\Big) \bigg ),
\end{split}
\end{equation*}
we rewrite the system \eqref{FVMunisch} in the form
   \begin{equation*}
      \mathcal F \big ( \bar u^{(n+1)}, \bar v^{(n+1)} \big ) + \mathcal G \big ( \bar u^{(n+1)}, \bar v^{(n+1)} \big ) = y := \big ( \bar u^{(n)}, \bar v^{(n)} \big ).
   \end{equation*}
Moreover we see that setting $\mathcal O = B(0,r)\subset \mathbb R^2$ a ball centered at $(0,0)$ with a radius
   \begin{equation*}
      r > \sqrt{\big( U + \tfrac{\alpha}{\beta}V \big)^2 + \big( V + \tfrac{\beta}{\alpha}U \big)^2}
   \end{equation*}
we fulfill all the assumptions of Theorem \cite[Theorem 3.1, page 16]{deimling}.\\
For the continuous function $\mathcal H:[0,1]\times B(0,r)\rightarrow \mathbb R^2$ given by
   \begin{equation*}
      \mathcal H(\lambda,\bar u^{(n+1)}, \bar v^{(n+1)}) = \mathcal F(\bar u^{(n+1)}, \bar v^{(n+1)}) + \lambda \mathcal G(\bar u^{(n+1)}, \bar v^{(n+1)}),
   \end{equation*}
we see that
   \begin{equation*}
      d \big (\mathcal F + \lambda \mathcal G, B, (\bar u^{(n)}, \bar v^{(n)}) \big ) = d \big (\mathcal H(\lambda), B, (\bar u^{(n)}, \bar v^{(n)}) \big )
   \end{equation*}
for all $\lambda \in [0,1]$ and $(\overline u^{(n)}),\,(\overline v^{(n)})$ such that $0 \leqslant \bar u^{(n)} \leqslant U - \tfrac{\alpha}{\beta}V$ and $0 \leqslant \bar v^{(n)} \leqslant V - \tfrac{\beta}{\alpha}U$ for $n \in \on$. On the other hand we deduce from \cite[Theorem 3.1, page 16]{deimling} (d1) that
\begin{equation}\label{FVMdeg_id}
   d \big (\mathcal H(0), B, (\bar u^{(n)}, \bar v^{(n)}) \big ) = 1.
\end{equation}
In view of \cite[Theorem 3.1, page 16]{deimling} (d3) and (d4), \eqref{FVMdeg_id} implies that the equality
\begin{equation*}
   \mathcal F(\bar u^{(n+1)}, \bar v^{(n+1)}) + \mathcal G(\bar u^{(n+1)}, \bar v^{(n+1)}) = \big ( \bar u^{(n)}, \bar v^{(n)} \big)
\end{equation*}
has a solution or, in other words, that there exists a solution of \eqref{FVMunisch}. The uniqueness of this solution immediately follows from Proposition \ref{FVMcomprinciple}.

\vskip 10pt
We can prove in the same way the existence and uniqueness of the solution of the system \eqref{FVMscheme}. Indeed, we rewrite \eqref{FVMscheme} in the form
   \begin{equation}\label{FVMeqtilde}
      \widetilde{\cal F}\Big( (u^{(n+1)}_K)_{K\in\cal M}, (v^{(n+1)}_K)_{K\in\cal M} \Big) + \widetilde{\cal G}\Big( (u^{(n+1)}_K)_{K\in\cal M}, (v^{(n+1)}_K)_{K\in\cal M} \Big) =\ \Big( (u^{(n)}_K)_{K\in\cal M}, (v^{(n)}_K)_{K\in\cal M} \Big),
   \end{equation}
where $\widetilde{\cal F}, \widetilde{\cal G}: \mathbb R^{2\Theta}\rightarrow \mathbb R^{2\Theta}$, with $\Theta$ the number of control volumes for the discretization $\cal D$, are continuous functions given by
\begin{equation*}
\begin{split}
   & \widetilde{\cal F}\Big( (u^{(n)}_K)_{K \in \cal M}, (v^{(n)}_K)_{K \in \cal M} \Big) = \Big( (u^{(n)}_K)_{K \in \cal M}, (v^{(n)}_K)_{K \in \cal M} \Big),\\
   & \widetilde{\cal G}\Big( (u^{(n)}_K)_{K \in \cal M}, (v^{(n)}_K)_{K \in \cal M} \Big) = \big( \mathcal W_1, \mathcal W_2 \big),
\end{split}
\end{equation*}
where
\begin{equation*}
   \rule{15pt}{0pt}\mathcal W_1 = -\frac{\dt^{(n)} a}{m_K} \sum_{L \in{\cal N}_K} T_{K|L} (u_L^{(n)} - u_K^{(n)}) + \dt^{(n-1)} \ha \big(r_A(u_K^{(n)}) - r_B(v_K^{(n)})\big),\rule{15pt}{0pt}
\end{equation*}
and where
\begin{equation*}
   \rule{15pt}{0pt} \mathcal W_2 = -\frac{\dt^{(n)} b}{m_K} \sum_{L \in{\cal N}_K} T_{K|L} (v_L^{(n)} - v_K^{(n)}) - \dt^{(n-1)} \hb \big(r_A(u_K^{(n)}) - r_B(v_K^{(n)})\big).\rule{15pt}{0pt}
\end{equation*}

We set $\widetilde{\mathcal O} = B(0,R)\subset \mathbb R^{2\Theta}$ a ball centered at zero with a radius
\begin{equation*}
   R > \sqrt{\Theta \big( U + \tfrac{\alpha}{\beta}V \big)^2 + \Theta \big( V + \tfrac{\beta}{\alpha}U \big)^2}.
\end{equation*}
Since $\Theta > 1$, we deduce from the discrete $L^\infty(Q_T)$ estimate of Theorem \ref{FVMlinf} that the equation \eqref{FVMeqtilde} does not have any solutions on $\partial \widetilde{\mathcal O}$. Applying again \cite[Theorem 3.1, page 16]{deimling} with $\widetilde{\cal H}(\lambda) = \widetilde{\cal F} + \lambda \widetilde{\cal G}$ and $\lambda\in[0,1]$ completes the proof of the following result.
\begin{theo}
   We suppose that the hypotheses $\mathcal H$ are satisfied. Let $\cal D$ be a discretization as in Definition \ref{FVMdiscretization}. Let  $(u_K^{(0)},v_K^{(0)})_{K\in{\cal M}}$ be given by \eqref{FVMschini}. Then there exists one and only one sequence
   \begin{equation*}
      (u_K^{(n+1)},v_K^{(n+1)})_{K\in{\cal M},\ n\in\on},
   \end{equation*}
   which satisfies \eqref{FVMscheme}, with the initial condition $(u_K^{(0)},v_K^{(0)})_{K\in{\cal M}}$.\hskip 5pt $\blacksquare$
\end{theo}

\section{Convergence proof with $k$ fixed}\label{FVMfixed_k}

We begin with the discrete version of $L^2(Q_T)$ estimates of the gradient of the approximate solutions.
\begin{prop}\label{FVMfixkgrad}
   We suppose that the hypotheses $\mathcal H$ are satisfied. Let $\cal D$ be a discretization as in Definition \ref{FVMdiscretization}. Let \eqref{FVMschini} and \eqref{FVMscheme} give the sequences $(u_K^{(0)},v_K^{(0)})_{K\in{\cal M}}$ and $(u_K^{(n+1)},v_K^{(n+1)})_{K\in{\cal M},\ n\in\on}$, respectively. Then, there exists a constant $C_k > 0$, which does not depend on $\cal D$, but which depend on all the data of the continuous \Pp\ (namely, the constants $\alpha,\ \beta,\ U,\ V$ including $k$ and the functions $r_A,\ r_B$), such that
   \begin{multline}\label{FVMfixkgradu}
      \rule{20pt}{0pt} \frac{1}{2} \sum_{K \in {\cal M}} m_K \Big( u_K^{(N+1)} \Big)^2 - \frac{1}{2} \sum_{K \in {\cal M}} m_K \Big( u_K^{(0)} \Big)^2 + a \sum_{n=0}^N \dt^{(n)} \sum_{(K,L)\in{\cal E}} T_{K|L} \Big(u_L^{(n+1)} - u_K^{(n+1)}\Big)^2 \leqslant C_k \rule{20pt}{0pt}
   \end{multline}
   and
   \begin{multline}\label{FVMfixkgradv}
      \rule{20pt}{0pt} \frac{1}{2} \sum_{K \in {\cal M}} m_K \Big( v_K^{(N+1)} \Big)^2 - \frac{1}{2} \sum_{K \in {\cal M}} m_K \Big( v_K^{(0)} \Big)^2 + b \sum_{n=0}^N \dt^{(n)} \sum_{(K,L)\in{\cal E}} T_{K|L} \Big(v_L^{(n+1)} - v_K^{(n+1)}\Big)^2 \leqslant C_k \rule{20pt}{0pt}
   \end{multline}
\end{prop}
{\hskip -\parindent \bf Proof \hskip 5pt}
   For the sake of simplicity we only present the proof for the $u$-component. We multiply the first equation in the finite volume scheme \eqref{FVMscheme} by $u_K^{(n+1)}$ and sum the result over all $K \in \cal M$ and over all $n\in \on$ to obtain
   \begin{equation}\label{FVM1gwiazdka}
      \mathcal S_1 + \mathcal S_2 + \mathcal S_3 = 0,
   \end{equation}
   where
   \begin{gather*}
      \mathcal S_1 = \sum_{n=0}^N \sum_{K \in {\cal M}} m_K \Big( u_K^{(n+1)} -  u_K^{(n)}\Big) u_K^{(n+1)},\\
      \mathcal S_2 = - a \sum_{n=0}^N \dt^{(n)} \sum_{K \in {\cal M}} \sum_{L \in{\cal N}_K} T_{K|L} \Big(u_L^{(n+1)} - u_K^{(n+1)}\Big) u_K^{(n+1)}, \\
      \mathcal S_3 = \hat \alpha \sum_{n=0}^N \dt^{(n)} \sum_{K \in {\cal M}} \Big( r_A \big( u_K^{(n+1)} \big) - r_B \big( v_K^{(n+1)} \big) \Big) u_K^{(n+1)}.
   \end{gather*} \medskip

   \hskip -\parindent Since
   \begin{equation*}
      \Big( u_K^{(n+1)}\Big)^2 - u_K^{(n)} u_K^{(n+1)} = \frac{1}{2} \Big( u_K^{(n+1)}\Big)^2 - \frac{1}{2} \Big( u_K^{(n)}\Big)^2 + \frac{1}{2} \Big( u_K^{(n+1)} - u_K^{(n)} \Big)^2,
   \end{equation*}
   we deduce that
   \begin{multline}\label{FVMS1}
      \mathcal S_1 = \frac{1}{2} \sum_{n=0}^N \sum_{K \in {\cal M}} m_K \bigg( \Big( u_K^{(n+1)} \Big)^2 - \Big( u_K^{(n)} \Big)^2 \bigg) + \frac{1}{2} \sum_{n=0}^N \sum_{K \in {\cal M}} m_K \Big( u_K^{(n+1)} - u_K^{(n)} \Big)^2 \\
      \geqslant \frac{1}{2} \sum_{n=0}^N \sum_{K \in {\cal M}} m_K \bigg( \Big( u_K^{(n+1)} \Big)^2 - \Big( u_K^{(n)} \Big)^2 \bigg),
   \end{multline}
   All the terms in the sum on $n$ on the right hand side of \eqref{FVMS1} simplify except for the first and the last ones. We have that
   \begin{equation}\label{FVM2gwiazdka}
      \mathcal S_1 \geqslant \frac{1}{2} \sum_{K \in {\cal M}} m_K \Big( u_K^{(N+1)} \Big)^2 - \frac{1}{2} \sum_{K \in {\cal M}} m_K \Big( u_K^{(0)} \Big)^2.
   \end{equation}
   \medskip

   \hskip -\parindent We can perform a discrete integration by parts to obtain
   \begin{equation}\label{FVM3gwiazdka}
      \mathcal S_2 = a \sum_{n=0}^N \dt^{(n)} \sum_{(K,L)\in{\cal E}} T_{K|L} \Big(u_L^{(n+1)} - u_K^{(n+1)}\Big)^2.
   \end{equation}\medskip

   \hskip -\parindent Finally we use the hypothesis $\mathcal H\, 4$ and the inequalities in \eqref{FVMeqlinf} to estimate the last term, namely
   \begin{multline}\label{FVM4gwiazdka}
      - \mathcal S_3 \leqslant \hat \alpha \sum_{n=0}^N \dt^{(n)} \sum_{K \in {\cal M}} m_K \Big( r_A \big( U + \frac{\alpha}{\beta}V \big) + r_B \big( V + \frac{\beta}{\alpha} U \big) \Big)\big( U + \frac{\alpha}{\beta}V \big) \\
      \leqslant \alpha k C \sum_{n=0}^N \dt^{(n)} \sum_{K \in {\cal M}} m_K = \alpha k C |\Omega| T,
   \end{multline}
   with some positive constant $C$.\medskip

   \hskip -\parindent Identities \eqref{FVM1gwiazdka} and \eqref{FVM3gwiazdka} together with the inequalities \eqref{FVM2gwiazdka} and \eqref{FVM4gwiazdka} immediately give \eqref{FVMfixkgradu}. Since the argument in the case of the $v$-component is similar, we omit the proof. {\hskip 5pt $\blacksquare$}

\subsection{Space and time translates of approximate solutions} \label{FVMfixksection}

We now turn to the space translates estimates. We use here methods which have been presented for example by Eymard, Gutnic and Hilhorst \cite{eymard_gutnic} and by Eymard, Gallou{\"e}t, Hilhorst and Slimane \cite{eymard_slimane}. The results of the current and the next subsection together with the technical Proposition \ref{FVMconskol} will imply the relative compactness of the sequence of approximate solutions.
\begin{prop}[Space translates estimates]\label{FVMfixkspace}
   We assume that
   \begin{enumerate}
      \item ${\cal D} = \Big({\cal M}, {\cal E}, (x_K)_{K \in {\cal T}},  (t^{(n)})_{n\in\onp}\Big)$ is an admissible discretization of $Q_T$ in the sense of Definition \ref{FVMdiscretization},
      \item the hypotheses $\mathcal H$ and assumptions \eqref{FVMhypsup} are satisfied,
      \item the functions $(u_{\cal D}$ and $v_{\cal D})$ are derived from the scheme \eqref{FVMschini}\,--\,\eqref{FVMscheme} and given by the formulas \eqref{FVMuapp}.
   \end{enumerate}
   Then there exists a positive constant $C_k$, which does not depend on $\cal D$, but depends on all the data of the continuous \Pp,\ including $k$, such that
   \begin{equation}
      \rule{15pt}{0pt} \int_0^T \int_{\Omega_\xi} \big(u_{\cal D}(x+\xi,t) - u_{\cal D}(x,t)\big)^2 \dx \dy{t} \leqslant C_k |\xi| \big(2\, \size({\cal D} \big) +|\xi|),\rule{15pt}{0pt} \label{FVMfixkspaceu}
   \end{equation}
   and
   \begin{equation}
      \rule{15pt}{0pt}\int_0^T \int_{\Omega_\xi} \big(v_{\cal D}(x+\xi,t) - v_{\cal D}(x,t)\big)^2 \dx \dy{t} \leqslant C_k |\xi| \big(2\, \size({\cal D}) \big) +|\xi|),\rule{15pt}{0pt}\label{FVMfixkspacev}
   \end{equation}
   for all $\xi\in\mathbb R^d$ and for $\Omega_\xi$ defined as in Proposition \ref{FVMconskol}.
\end{prop}
{\hskip -\parindent \bf Proof \hskip 5pt}
   Inequalities \eqref{FVMfixkspaceu} and \eqref{FVMfixkspacev} follow from the estimates \eqref{FVMfixkgradu} and \eqref{FVMfixkgradv}, respectively. We refer to \cite[Lemma 3.3]{ciarlet_lions_VII} for a details.
   {\hskip 5pt $\blacksquare$}

\begin{prop}[Time translates estimates]\label{FVMfixktime}
   Let the assumptions of Proposition \ref{FVMfixkspace} are satisfied. Then, there exists some constant $C_k > 0$, which does not depend on $\cal D$, but which depend on all the data including $k$, such that
   \begin{equation}\label{FVMfixktimeu}
      \rule{20pt}{0pt} \int_{\Omega \times (0,T-\tau)} \big(u_{\cal D}(x,t+\tau) - u_{\cal D}(x,t) \big)^2 \dx \dy{t} \leqslant C_k \big(\size({\cal D}) +\tau \big) \rule{20pt}{0pt}
   \end{equation}
   and
   \begin{equation}\label{FVMfixktimev}
      \rule{20pt}{0pt} \int_{\Omega \times (0,T-\tau)} \big(v_{\cal D}(x,t+\tau) - v_{\cal D}(x,t) \big)^2 \dx \dy{t} \leqslant C_k \big(\size({\cal D}) +\tau \big), \rule{20pt}{0pt}
   \end{equation}
   for all $\tau\in(0,T)$.
 \end{prop}
{\hskip -\parindent \bf Proof \hskip 5pt}
   In order to apply Lemma \ref{FVMestttt} (see appendix), we follow the same steps as in \cite[Lemma 5.5]{eymard_gutnic}. The only difference appear in the nonlinear part of the equations. However, these can be easily estimated using the regularity properties of functions $r_A(\cdot)$ and $r_B(\cdot)$, as well as $L^\infty$ estimates \eqref{FVMeqlinf} in Theorem \ref{FVMlinf}.\hskip 5pt $\blacksquare$

\subsection{Convergence proof}

In this section, we state convergence results with $k$ fixed. This differs from next section where we will introduce additional hypotheses about the nonlinear reaction terms and obtain convergence results which permit us to pass to the limit as $k\rightarrow\infty$.

\begin{theo}\label{FVMfixkconvth}
   We suppose that the hypotheses $\mathcal H$ are satisfied. Let $(u_{{\cal D}},v_{{\cal D}})$ be the approximate  solution defined by \eqref{FVMschini}, \eqref{FVMscheme} and \eqref{FVMuapp}. There exist a pair of functions $(u^k,v^k)$ and a sequence $(u_{{\cal D}_m},v_{{\cal D}_m})_{m\in\N}$ of $(u_{{\cal D}},v_{{\cal D}})$ such that
  \begin{equation*}
   (u_{{\cal D}_m},v_{{\cal D}_m})_{m\in\N} \mbox{~~converges to~~} (u^k,v^k)\in \big(L^2\big(0,T; H^1(\Omega) \big) \big)^2
   \end{equation*}
   strongly in $L^2\big(Q_T\big)$ as  $\size({\cal D}_m)$ tends to zero.
   The function pair $(u^k,v^k)$ is a weak solution of \Pp\ in the sense of Definition \ref{FVMweakdefi}.
\end{theo}\medskip

\noindent Since \Pp\ is a uniformly parabolic system, $(u^k,v^k)$  must coincide with the unique classical solution
of  \Pp\ . This immediately yields the following result.
\begin{coro}
The pair $(u_{{\cal D}},v_{{\cal D}})$ converges to the unique classical solution $(u^k,v^k)$ of \Pp\ as  $\size({\cal D})$ tends to zero. \hskip 5pt $\blacksquare$
\end{coro}\medskip

{\hskip -\parindent \bf Proof of Theorem \ref{FVMfixkconvth} \hskip 5pt}
   In view of the estimates \eqref{FVMfixkspaceu}, \eqref{FVMfixktimeu} and Pro\-po\-sition \ref{FVMconskol} which is a consequence of the Fr\'echet--Kolmogorov The\-o\-rem \cite[Theorem IV.25, page 72]{brezis}, we deduce the relative compactness of the set $\big( u_{{\cal D}}\big)$ so that there exists a sequence of $\big( u_{{\cal D}_m}\big)_{m=1}^\infty$ and a  function ${\cal U}_k$, such that $u_{{\cal D}_m}\rightarrow {\cal U}_k$ strongly in $L^2(Q_T)$ and weakly in $L^2\big(0,T;H^1(\Omega)\big)$, as $m\rightarrow\infty$.

   The same conclusion holds for the $v$-component. Indeed, the ine\-qua\-lities \eqref{FVMfixkspacev} and \eqref{FVMfixktimev} permit to apply the compactness result in Proposition \ref{FVMconskol} for the sequence $\big( v_{{\cal D}_m}\big)_{m=1}^\infty$. There exists a function ${\cal V}_k$ such that $v_{{\cal D}_m}\rightarrow {\cal V}_k$ strongly in $L^2(Q_T)$ and weakly in $L^2\big(0,T;H^1(\Omega)\big)$ as $m\rightarrow\infty$.

   Next we show that $(\mathcal U_k$, $\mathcal V_k)$ is a weak solution of \Pp, in the sense of Definition \ref{FVMweakdefi}. Since the proof for the $v$-component is similar, we only present here the detailed proof in case of the $u$-component.\\
   Let $\psi \in \Psi$, where $\Psi$ is the class of test functions from Definition \ref{FVMweakdefi}.
   We multiply the first equation of \eqref{FVMscheme} by $\psi(x_K,t^{(n)})$, where $\psi \in \Psi$. Then we sum over all $K\in\cal M$ and $n\in\{0\ldots N-1\}$ to obtain
   \begin{equation*}
      {\cal T}_{1m}^u - {\cal T}_{2m}^u + {\cal T}_{3m}^u = 0,
   \end{equation*}
   where
   \begin{align*}
      {\cal T}_{1m}^u & = \sum_{n=0}^{N-1} \sum_{K\in {\cal M}} m_K \big( u_K^{(n+1)} - u_K^{(n)} \big) \psi(x_K,t^{(n)}),\\
      {\cal T}_{2m}^u & = a \sum_{n=0}^{N-1} \dt^{(n)} \sum_{K\in {\cal M}} \sum_{L\in{\cal N}_K} T_{K|L} \big(u_L^{(n+1)} - u_K^{(n+1)} \big) \psi(x_K,t^{(n)}),\\
      {\cal T}_{3m}^u & = \sum_{n=0}^{N-1} \dt^{(n)} \sum_{K\in {\cal M}} m_K \ha \Big(r_A(u_K^{(n+1)}) - r_B(v_K^{(n+1)}) \Big)\psi(x_K,t^{(n)}).
   \end{align*}

   The complete proof, that
   \begin{equation*}
      \lim_{m \rightarrow \infty} {\cal T}_{1m}^u = -\int_{\Omega} u_0(x)\psi(x,0) \dx - \int_0^T \int_\Omega {\cal U}_k(x,t) \psi_t(x,t) \dx \dy{t}
   \end{equation*}
   and
   \begin{equation*}
      \lim_{m \rightarrow \infty} {\cal T}_{2m}^u = - a \int_0^T \int_\Omega {\cal U}_k(x,t) \Delta \psi(x,t) \dx \dy{t},
   \end{equation*}
   can by found in \cite[Lemma 5.5]{eymard_gutnic}. Let us focus on the proof that
   \begin{equation*}
      \lim_{m \rightarrow \infty} {\cal T}_{3m}^u = \alpha k \int_0^T \int_\Omega \big( r_A({\cal U}_k) - r_B({\cal V}_k) \big) \psi(x,t) \dx \dy{t}.
   \end{equation*}
   We write
   \begin{align*}
      \ha \sum_{n=0}^{N-1} \dt^{(n)} & \sum_{K\in {\cal M}} m_K \Big(r_A(u_K^{(n+1)}) - r_B(v_K^{(n+1)}) \Big)\psi(x_K,t^{(n)}) \\
      & - \ha \int_0^T \int_\Omega \psi(x,t) \big( r_A({\cal U}_k) - r_B({\cal V}_k) \big) \dx \dy{t}\\
      =&\, \ha \sum_{n=0}^{N-1} \sum_{K\in {\cal M}} \int_{t^{(n)}}^{t^{(n+1)}} \int_K \Big(r_A(u_K^{(n+1)}) - r_B(v_K^{(n+1)}) \Big)\psi(x_K,t^{(n)}) \dx \dy{t}\\
      &- \ha \sum_{n=0}^{N-1} \sum_{K\in {\cal M}} \int_{t^{(n)}}^{t^{(n+1)}} \int_K \psi(x,t) \big( r_A({\cal U}_k) - r_B({\cal V}_k) \big) \dx \dy{t} \\
      &- \ha \sum_{K\in {\cal M}} \int_{t^{(N)}}^{t^{(N+1)}} \int_K \psi(x,t) \big( r_A({\cal U}_k) - r_B({\cal V}_k) \big) \dx \dy{t}.
   \end{align*}
   Thanks to the regularity of the function $\psi$, the last sum above converges to zero . Moreover,
   \begin{equation}\label{FVMfixth4}
   \begin{split}
      \ha \sum_{n=0}^{N-1} &\sum_{K\in {\cal M}} \int_{t^{(n)}}^{t^{(n+1)}} \int_K \Big(r_A(u_K^{(n+1)}) - r_B(v_K^{(n+1)}) \Big)\psi(x_K,t^{(n)}) \dx \dy{t}\\
      &- \ha \sum_{n=0}^{N-1} \sum_{K\in {\cal M}} \int_{t^{(n)}}^{t^{(n+1)}} \int_K \psi(x,t) \big( r_A(\mathcal U_k) - r_B(\mathcal V_k) \big) \dx \dy{t} \\
      =&\, \ha \sum_{n=0}^{N-1} \sum_{K\in {\cal M}} \int_{t^{(n)}}^{t^{(n+1)}} \int_K \big( \psi(x_K,t^{(n)}) - \psi(x,t) \big) \times \\
      &\rule{120pt}{0pt}\times (r_A(u_K^{(n+1)}) - r_B(v_K^{(n+1)}) \dx \dy{t}\\
      &+ \ha \sum_{n=0}^{N-1} \sum_{K\in {\cal M}} \int_{t^{(n)}}^{t^{(n+1)}} \int_K \psi(x,t) \big(r_A(u_K^{(n+1)}) - r_A({\cal U}_k)\big) \dx \dy{t}\\
      &- \ha \sum_{n=0}^{N-1} \sum_{K\in {\cal M}} \int_{t^{(n)}}^{t^{(n+1)}} \int_K \psi(x,t) \big(r_B(v_K^{(n+1)}) - r_B({\cal V}_k)\big) \dx \dy{t}.
   \end{split}
   \end{equation}
   Next we show that the three terms above tend to zero as $m\rightarrow \infty$. First we take their absolute value and apply the triangle inequality. The Cauchy-Schwarz inequality applied to the first sum of the right hand side of \eqref{FVMfixth4} yields
   \begin{gather*}
      \ha \sum_{n=0}^{N-1} \sum_{K\in {\cal M}} \int_{t^{(n)}}^{t^{(n+1)}} \int_K \big| \psi(x_K,t^{(n)}) - \psi(x,t) \big| \big|r_A(u_K^{(n+1)}) - r_B(v_K^{(n+1)}\big| \dx \dy{t}\\
      \leqslant \bigg( \sum_{n=0}^{N-1} \sum_{K\in {\cal M}} \int_{t^{(n)}}^{t^{(n+1)}} \int_K \big( \psi(x_K,t^{(n)}) - \psi(x,t) \big)^2 \dx \dy{t} \bigg)^{1/2} \times \\
      \times \hat \alpha \bigg( \sum_{n=0}^{N-1} \sum_{K\in {\cal M}} \int_{t^{(n)}}^{t^{(n+1)}} \int_K \big(r_A(u_K^{(n+1)}) - r_B(v_K^{(n+1)}\big)^2 \dx \dy{t} \bigg)^{1/2} .
   \end{gather*}
   The first term of above product converges to zero, as $m\rightarrow \infty$, since $\psi(x,t)$ is smooth enough. The second term is bounded. Indeed, it is sufficient to remark that $r_A(u_K^{(n+1)})$, $r_B(v_K^{(n+1)})$ are bounded for all $K \in \mathcal M$ and $n \in \on$. The last two terms in \eqref{FVMfixth4} are similar and we show how the proof goes with the first one. Indeed,
   \begin{gather*}
      \ha \sum_{n=0}^{N-1} \sum_{K\in {\cal M}} \int_{t^{(n)}}^{t^{(n+1)}} \int_K |\psi(x,t)| \big| r_A(u_K^{(n+1)}) - r_A({\cal U}_k) \big| \dx \dy{t}\\
      \leqslant \, \ha \|\psi \|_{L^\infty(Q_T)} \sum_{n=0}^{N-1} \sum_{K\in {\cal M}} \int_{t^{(n)}}^{t^{(n+1)}} \int_K  \big| r_A(u_K^{(n+1)}) - r_A({\cal U}_k) \big| \dx \dy{t}\\
      \leqslant \,\ha \|\psi \|_{L^\infty(Q_T)} \sum_{n=0}^{N-1} \sum_{K\in {\cal M}} \int_{t^{(n)}}^{t^{(n+1)}} \int_K  \| r_A' \|_{L^\infty(Q_T)} \big| u(x_K,t^{(n)}) - {\mathcal U_k}(x,t) \big| \dx \dy{t}\\
      \leqslant \, \ha \|\psi \|_{L^\infty(Q_T)} \bigg( \sum_{n=0}^{N-1} \sum_{K\in {\cal M}} \int_{t^{(n)}}^{t^{(n+1)}} \int_K \| r_A' \|_{L^\infty(Q_T)}^2 \dx \dy{t}\bigg)^{1/2}\ \times \\
      \times \bigg( \sum_{n=0}^{N-1} \sum_{K\in {\cal M}} \int_{t^{(n)}}^{t^{(n+1)}} \int_K \big| u(x_K,t^{(n)}) - {\mathcal U_k}(x,t) \big|^2 \dx \dy{t} \bigg)^{1/2}.
   \end{gather*}
   The $L^{\infty}$ norm of the function $r_A(x)$ is taken over the finite interval $[0,U+\frac{\alpha}{\beta}V]$. Because $r_A(x)$ is of class $C^1(\mathcal R)$ the first term of the above product is bounded. The second one converges to zero since $u_{{\cal D}_m} \rightarrow {\cal U}_k$ as $m\rightarrow \infty$ in $L^2(Q_T)$.{\hskip 5pt $\blacksquare$}

\section{The case that $k$ tends to infinity}\label{FVMunifcase}

In order to prove the convergence of the finite volume scheme when $size~(\cal D)$ tends to zero and $k$ tends to infinity, we impose some additional conditions on the nonlinear terms $r_A(x)$ and $r_B(x)$.
At first we prove a counterpart of Proposition \ref{FVMfixkgrad}.
\begin{prop}\label{FVMestigrad}
   Let us assume hypotheses $\mathcal H$. Moreover, we assume that the functions $r_A(x),\ r_B(x)$ satisfy
   \begin{equation}\label{FVMhypsup}
   \begin{split}
      & r_\kappa\in C^1(\mathbb R),\ r_\kappa'(\cdot) > 0 \text{\quad on \quad}(0,+\infty),\\
      & r_\kappa(0) = 0, \text{\quad and \quad} \limsup_{s\rightarrow 0^+} \frac{s r_\kappa'(s)}{r_\kappa(s)} < \infty,
   \end{split}
   \end{equation}
   where $\kappa\in\{A,\ B\}$. Let ${\cal D} = \Big({\cal M}, {\cal E}, (x_K)_{K \in {\cal T}}, (t^{(n)})_{n\in\onp}\Big)$ be an admissible discretization of $Q_T$ in the sense of Definition \ref{FVMdiscretization}, and the sequences $(u_K^{(0)},v_K^{(0)})_{K\in{\cal M}}$ and $\big(u_K^{(n+1)},v_K^{(n+1)}\big)_{K\in{\cal M},\ n\in\on}$ are given by \eqref{FVMschini} and \eqref{FVMscheme}, respectively. Then, there exists some positive constant $C$ which is independent of the discretization $\cal D$ and of the reaction rate $k$, such that
   \begin{equation}\label{FVMeql2h1d}
      \sum_{n=0}^N \dt^{(n)} \bigg( \sum_{(K,L)\in{\cal E}} T_{K|L} \big( u_L^{(n+1)} - u_K^{(n+1)}\big)^2 + \sum_{(K,L)\in{\cal E}} T_{K|L} \big(v_L^{(n+1)} - v_K^{(n+1)}\big)^2\bigg) \leqslant C
   \end{equation}
   and
   \begin{equation}\label{FVMeqchose}
      k \sum_{n=0}^N \dt^{(n)} \sum_{K\in{\cal M}}  m_K \big(r_A(u_K^{(n+1)}) - r_B(v_K^{(n+1)})\big)^2 \leqslant C.
   \end{equation}
\end{prop}

\begin{rema}
   Observe that the condition \eqref{FVMhypsup} holds, for example, in the case that the rate functions $r_A(s)$ and $r_B(s)$ behave like $s^\gamma$, for some positive $\gamma$, whenever $s \rightarrow 0^+$.
\end{rema}

{\hskip -\parindent \bf Proof of Proposition \ref{FVMestigrad} \hskip 5pt}
   Let $(a,b)\in (\mathbb R^+)^2$ be such that $r_A(a) = r_B(b)$. We define two functions
   \begin{equation*}
   \begin{split}
      V_A(s) = \frac{1}{\alpha} \Bigg( s\ln \frac{r_A(s)}{r_A(a)} + \int_s^a \frac{\sigma r_A'(\sigma)}{r_A(\sigma)} \dy{\sigma} \Bigg), \\
      V_B(s) = \frac{1}{\beta}  \Bigg( s\ln \frac{r_B(s)}{r_B(b)} + \int_s^b \frac{\sigma r_B'(\sigma)}{r_B(\sigma)} \dy{\sigma} \Bigg),
   \end{split}
   \end{equation*}
   which are continuous on $\mathbb R^+$ because of hypotheses $\mathcal H$ and the assumptions \eqref{FVMhypsup}. We can extend these functions to also be continuous at $s=0$. To do so for the function $V_A(s)$ we check that the hypotheses \eqref{FVMhypsup} give the integrability of $\displaystyle \frac{s r_A'(s)}{r_A(s)}$ on the interval $[0,a]$ and we pass to the limit
   \begin{equation*}
      \lim_{s\rightarrow 0^+} V_A(s) \stackrel{H}{=} \frac{1}{\alpha} \int_0^a \frac{\sigma r_A'(\sigma)}{r_A(\sigma)} \dy{\sigma} - \frac{1}{\alpha} \lim_{s\rightarrow 0^+} s\cdot \frac{s\,r_A'(s)}{r_A(s)} < \infty,
   \end{equation*}
   where we have applied de l'Hospital theorem as formulated in \cite[Theorem 2, p.\ 174]{nikolsky1}.

   For a given $\eps\in(0,1)$ and $n\in\on$ we consider
   \begin{equation*}
      A^{(n+1)}(\eps) = \sum_{K\in{\cal M}} m_K \Big( V_A\big( u_K^{(n+1)}+\eps\big) - V_A\big(
      u_K^{(n)}+\eps\big)\Big).
   \end{equation*}
   Since $\displaystyle V_A''(s) = \frac {r_A'(s)} {\alpha\,r_A(s)} >0$ for all $s \in \mathbb R^+$ we deduce that
   \begin{gather*}
      V_A \big( u_K^{(n)} + \eps \big) - V_A \big(u_K^{(n+1)} + \eps \big) = \big(u_K^{(n)} - u_K^{(n+1)}\big) V_A' \big(u_K^{(n+1)} + \eps \big)+\frac{1}{2} \big(u_K^{(n)} - u_K^{(n+1)}\big)^2 V_A'' \big(s \big)\\
         \geqslant \big(u_K^{(n)} - u_K^{(n+1)}\big)V_A' \big(u_K^{(n+1)} + \eps \big)\hskip 5pt.
   \end{gather*}
   As a consequence
   \begin{equation*}
      A^{(n+1)}(\eps) \leqslant \sum_{K\in{\cal M}} m_K
      \big( u_K^{(n+1)}-  u_K^{(n)}\big) V_A'\big( u_K^{(n+1)}+\eps \big).
   \end{equation*}
   We substitute $m_K \big( u_K^{(n+1)}-  u_K^{(n)}\big)$ from the scheme \eqref{FVMscheme}, which yields
   \begin{equation*}
      A^{(n+1)}(\eps) \leqslant  A_1^{(n+1)}(\eps) + A_2^{(n+1)}(\eps),
   \end{equation*}
   with
   \begin{equation*}
      A_1^{(n+1)}(\eps) = - \dt^{(n)} \frac{a}{4} \sum_{(K,L)\in{\cal E}} T_{K|L}\big(u_L^{(n+1)} - u_K^{(n+1)}\big) \Big \{ V_A'\big( u_L^{(n+1)}+\eps\big) - V_A'\big( u_K^{(n+1)}+\eps \big)\Big \}
   \end{equation*}
   and
   \begin{equation*}
      A_2^{(n+1)}(\eps)= - \dt^{(n)}\sum_{K\in{\cal M}}  m_K \hat \alpha\, V_A'\big( u_K^{(n+1)}+\eps\big) \Big \{r_A(u_K^{(n+1)}) - r_B(v_K^{(n+1)})\Big \}.
   \end{equation*}
   Since there exists some constant $C >0$ such that $V_A''(s) \geqslant C$ for all $s\in [0,U+\frac\alpha \beta V]$ we can use the $L^\infty$ bound \eqref{FVMeqlinf} and the \textit{mean value theorem} to obtain
   \begin{equation*}
      A_1^{(n+1)}(\eps)\leqslant  - \dt^{(n)} C \frac{a}{4} \sum_{(K,L)\in{\cal E}} T_{K|L} \big(u_L^{(n+1)} - u_K^{(n+1)}\big)^2.
   \end{equation*}
   Following the same steps for the function
   \begin{equation*}
      B^{(n+1)}(\eps) = \sum_{K\in{\cal M}} m_K \Big( V_B\big( v_K^{(n+1)}+\eps\big) - V_B\big(v_K^{(n)}+\eps\big)\Big),
   \end{equation*}
   we arrive at
   \begin{equation}
      A^{(n+1)}(\eps)+B^{(n+1)}(\eps) \leqslant -\ C^{(n+1)} - D^{(n+1)}(\eps),
   \label{FVMapb}\end{equation}
   with
   \begin{equation*}
      C^{(n+1)} = \dt^{(n)} C \frac{a}{4} \sum_{(K,L)\in{\cal E}} T_{K|L} \big(u_L^{(n+1)} - u_K^{(n+1)} \big)^2 + \dt^{(n)} C \frac{b}{4} \sum_{(K,L)\in{\cal E}} T_{K|L} \big(v_L^{(n+1)} - v_K^{(n+1)} \big)^2,
   \end{equation*}
   and
   \begin{multline*}
      D^{(n+1)}(\eps) = \dt^{(n)} k \sum_{K\in{\cal M}} m_K \big(r_A(u_K^{(n+1)}) - r_B(v_K^{(n+1)})\big) \big(\ln r_A(u_K^{(n+1)}+\eps) - \ln r_B(v_K^{(n+1)}+\eps)\big),
   \end{multline*}
   where we used that $\displaystyle \alpha V_A'(s) = \ln \frac{r_A(s)}{r_A(a)}$ and that $\displaystyle \beta V_B'(s) = \ln \frac{r_B(s)}{r_B(b)}$.\\
   Let
   \begin{equation*}
      D_1^{(n+1)}(\eps) = \dt^{(n)} k \sum_{K\in{\cal M}}  m_K
      \big(r_A(u_K^{(n+1)}+\eps) - r_B(v_K^{(n+1)}+\eps)\big) \big(\ln r_A(u_K^{(n+1)}+\eps) -
      \ln r_B(v_K^{(n+1)}+\eps)\big).
   \end{equation*}
   Since the inequality
   \begin{equation*}
      (c-d)(\ln c - \ln d)\geqslant \frac{(c-d)^2}{c+d}
   \end{equation*}
   holds for all $\displaystyle (c,d)\in \Big(0,U+\frac \alpha \beta V +1 \Big) \times \Big(0,V+\frac \beta \alpha U +1\Big)$, then
   \begin{equation*}
      D_1^{(n+1)}(\eps) \geqslant \dt^{(n)} \frac{k}{C_b} \!\sum_{K\in{\cal M}}\! m_K \big(r_A(u_K^{(n+1)}+\eps) - r_B(v_K^{(n+1)}+\eps)\big)^2,
   \end{equation*}
   where $C_b$ is an upper bound for $r_A(c) + r_B(d)$ with $\displaystyle (c,d)\in \big(0,U+\frac \alpha \beta V +1\big)\times \big(0,V+\frac \beta \alpha U +1\big)$. Such bounds exist in view of Theorem \ref{FVMlinf} and regularity of the functions $r_A(x)$ and $r_B(x)$. Let us define
   \begin{equation*}
      E^{(n+1)}(\eps) := \dt^{(n)} \frac{k}{C_b} \!\sum_{K \in {\cal M}} \! m_K \big(r_A(u_K^{(n+1)} + \eps) - r_B(v_K^{(n+1)} + \eps)\big)^2.
   \end{equation*}
   The assumptions \eqref{FVMhypsup} and Lemma \ref{FVMoszacLN} in the appendix, imply that for all $K \in \cal M$ and $\eps >0$ small enough, there exist constants $C_1,\ C_2$ such that
   \begin{equation*}
      \big|r_A(u_K^{(n+1)}+\eps)-r_A(u_K^{(n+1)})\big|\leqslant C_1 \eps,\quad \big|r_B(v_K^{(n+1)}+\eps)-r_B(v_K^{(n+1)})\big|\leqslant C_1 \eps,
   \end{equation*}
   and
   \begin{equation*}
      \big|\ln r_A(u_K^{(n+1)} + \eps) \big| \leqslant C_2 \big( |\ln \eps|+1 \big),\quad \big|\ln r_B(v_K^{(n+1)} + \eps) \big| \leqslant C_2 \big( |\ln \eps|+1 \big).
   \end{equation*}
   Then
   \begin{equation*}
      - D^{(n+1)} \leqslant C \eps |\Omega| k \big( |\ln \eps| + 1 \big) - D^{(n+1)}(\eps),
   \end{equation*}
   for some positive constant $C$. As a consequence
   \begin{align*}
      A^{(n+1)}(\eps) + B^{(n+1)}(\eps) &\leqslant -\, C^{(n+1)} - D^{(n+1)}(\eps)\\
      &\leqslant -\, C^{(n+1)} - E^{(n+1)}(\eps) + C \eps |\Omega| k \big( |\ln \eps| + 1 \big).
   \end{align*}
   Now it is possible to pass to the limit in \eqref{FVMapb}. We obtain
   \begin{equation*}
      A^{(n+1)}(0) + B^{(n+1)}(0) \leqslant - C^{(n+1)} - E^{(n+1)}(0).
   \end{equation*}
   which is
   \begin{gather*}
      \sum_{K \in {\cal M}} m_K \Big( V_A \big( u_K^{(n+1)}\big) - V_A \big( u_K^{(n)}\big) \Big) + \sum_{K \in {\cal M}} m_K \Big( V_B \big( v_K^{(n+1)}\big) - V_B \big( v_K^{(n)}\big) \Big) \\
      + \dt^{(n)} C \frac{a}{4} \sum_{(K,L)\in{\cal E}} T_{K|L} \big(u_L^{(n+1)} - u_K^{(n+1)} \big)^2 + \dt^{(n)} C \frac{b}{4} \sum_{(K,L)\in{\cal E}} T_{K|L} \big(v_L^{(n+1)} - v_K^{(n+1)} \big)^2 \\
      + \dt^{(n)} \frac{k}{C_b} \!\sum_{K \in {\cal M}} \! m_K \big(r_A(u_K^{(n+1)}) - r_B(v_K^{(n+1)})\big)^2 \leqslant 0.
   \end{gather*}
   Now we sum the above inequalities over $n\in\on$ to obtain
   \begin{gather*}
      \sum_{K\in{\cal M}} m_K \Big( V_A \big( u_K^{(N+1)}\big)+V_B \big( u_K^{(N+1)}\big)\Big)+ + C \frac{a}{4} \sum_{n=0}^N \dt^{(n)} \sum_{(K,L)\in{\cal E}} T_{K|L} \big(u_L^{(n+1)} - u_K^{(n+1)} \big)^2 \\
      + C \frac{b}{4} \sum_{n=0}^N \dt^{(n)} \sum_{(K,L)\in{\cal E}} T_{K|L} \big(v_L^{(n+1)} - v_K^{(n+1)} \big)^2 + \frac{k}{C_b} \sum_{n=0}^N \dt^{(n)} \!\sum_{K \in {\cal M}} \! m_K \big(r_A(u_K^{(n+1)}) - r_B(v_K^{(n+1)})\big)^2\\
      \leqslant \sum_{K\in{\cal M}} m_K \Big(V_A \big( u_K^{(0)}\big)+ V_B \big( u_K^{(0)}\big)\Big).
   \end{gather*}
   Since $V_A(s)$ and $V_B(s)$ are nonnegative and continuous on $[0,+\infty)$, and since $u_K^{(n+1)},\ v_K^{(n+1)}$ are nonnegative and bounded for all $K \in \mathcal M$ and $n\in\on$, this concludes the proof.{\hskip 5pt $\blacksquare$ \medskip}

\subsection{Space and time translates of the approximate solutions}

Since we have already presented the general methods in section \ref{FVMfixed_k}
we only give here some essential ideas, leaving out the details of the proofs.

We begin with a counterpart of Proposition \ref{FVMfixkspace}.
\begin{prop}[Space translates estimates]\label{FVMlempom1}
   Let us assume that
   \begin{enumerate}
      \item $\Big({\cal M}, {\cal E}, (x_K)_{K \in {\cal T}},  (t^{(n)})_{n\in\onp}\Big)$ is an admissible discretization of $Q_T$ in the sense of Definition \ref{FVMdiscretization},
      \item hypotheses $\mathcal H$ and assumptions \eqref{FVMhypsup} are satisfied,
      \item functions $(u_{\cal D}$ and $v_{\cal D})$ are derived from the scheme \eqref{FVMschini}\,--\,\eqref{FVMscheme} and given by the formulas \eqref{FVMuapp}.
   \end{enumerate}
   Then there exists a positive constant $C$ which is independent of $\cal D$ and $k$, such that
   \begin{equation}
      \rule{15pt}{0pt} \int_0^T \int_{\Omega_\xi} \big(u_{\cal D}(x+\xi,t) - u_{\cal D}(x,t)\big)^2 \dx \dy{t} \leqslant C |\xi| \big(2\, \size({\cal D} \big) +|\xi|),\rule{15pt}{0pt} \label{FVMspacetr1}
   \end{equation}
   and
   \begin{equation}
      \rule{15pt}{0pt}\int_0^T \int_{\Omega_\xi} \big(v_{\cal D}(x+\xi,t) - v_{\cal D}(x,t)\big)^2 \dx \dy{t} \leqslant C |\xi| \big(2\, \size({\cal D} \big) +|\xi|),\rule{15pt}{0pt}\label{FVMspacetr2}
   \end{equation}
   for all $\xi\in\mathbb R^d$ and for $\Omega_\xi$ defined as in Proposition \ref{FVMconskol}.
\end{prop}
{\hskip -\parindent \bf Proof \hskip 5pt}
   As it was in the proof of Proposition \ref{FVMfixkspace}, we refer to \cite{ciarlet_lions_VII}[Lemma 3.3] for a complete proof. The only difference is to apply the result of Proposition \ref{FVMestigrad} instead of that of Proposition \ref{FVMfixkgrad}. {\hskip 5pt $\blacksquare$ \medskip}

Let us now prove an analogue of Proposition \ref{FVMfixktime}
\begin{lemm}[Time translates estimate]\label{FVMlempom2}
   Let the assumptions \textit{1.\ 2.} and \textit{3.}\ of Lemma \ref{FVMlempom1} be satisfied. Set $\displaystyle w_{\cal D} = \frac 1 \alpha u_{\cal D}+ \frac 1 \beta v_{\cal D}$. Then there exists a positive constant $C$ which is independent of $\cal D$ and $k$, such that
   \begin{equation*}
      \int_{\Omega\times(0,T-\tau)} \big(w_{\cal D}(x,t+\tau)-w_{\cal D}(x,t)\big)^2 \dx \dy{t} \leqslant C \big(\size({\cal D}) +\tau \big),
   \end{equation*}
   for all $\tau\in(0,T)$.
\end{lemm}
{\hskip -\parindent \bf Proof \hskip 5pt}
   The proof is similar to that of Proposition \ref{FVMfixktime}. We present below the essential steps of the argument. \medskip \\
   We define
   \begin{equation*}
      \calA(t) := \int_\Omega \big(w_{\cal D}(x,t+\tau)- w_{\cal D}(x,t)\big)^2 \dx,
   \end{equation*}
   which can be easily transformed into
   \begin{align*}
      \calA(t) &= \sum_{K \in \mathcal M} m_K \big( w_K^{(n(t+\tau)+1)} - w_K^{(n(t)+1)} \big)^2 \\
         &= \sum_{k \in \mathcal M} \bigg( \big( w_K^{(n(t+\tau)+1)} - w_K^{(n(t)+1)} \big) \sum_{n=n(t)+1}^{n(t+\tau)} m_K \big(w_K^{(n+1)} - w_K^{(n)} \big) \bigg).
   \end{align*}
   Since
   \begin{equation*}
      w^{(n+1)}_K - w^{(n)}_K = \frac{1}{\alpha} \Big( u_K^{(n+1)} - u_K^{(n)} \Big) + \frac{1}{\beta} \Big( v_K^{(n+1)} - v_K^{(n)} \Big),
   \end{equation*}
   we can apply discrete integration by parts in the scheme \eqref{FVMscheme} to obtain
   \begin{gather*}
      \calA(t) = \frac{a}{\alpha} \sum_{n=n(t)+1}^{n(t+\tau)} \dt^{(n)} \sum_{(K,L) \in \cal E} T_{K|L} \big( u_L^{(n+1)} - u_K^{(n+1)} \big) \big( w_K^{(n(t+\tau)+1)} - w_L^{(n(t+\tau)+1)}\big)\\
      + \frac{a}{\alpha} \sum_{n=n(t)+1}^{n(t+\tau)} \dt^{(n)} \sum_{(K,L) \in \cal E} T_{K|L} \big( u_L^{(n+1)} - u_K^{(n+1)} \big) \big( w_L^{(n(t)+1)} - w_K^{(n(t)+1)}\big)\\
      \frac{b}{\beta} \sum_{n=n(t)+1}^{n(t+\tau)} \dt^{(n)} \sum_{(K,L) \in \cal E} T_{K|L} \big( v_L^{(n+1)} - v_K^{(n+1)} \big) \big( w_K^{(n(t+\tau)+1)} - w_L^{(n(t+\tau)+1)}\big)\\
      + \frac{b}{\beta} \sum_{n=n(t)+1}^{n(t+\tau)} \dt^{(n)} \sum_{(K,L) \in \cal E} T_{K|L} \big( v_L^{(n+1)} - v_K^{(n+1)} \big) \big( w_L^{(n(t)+1)} - w_K^{(n(t)+1)}\big).
   \end{gather*}
   Next we estimate the second term in the sum above, to obtain
   \begin{gather*}
      \frac{a}{\alpha} \sum_{n=n(t)+1}^{n(t+\tau)} \dt^{(n)} \sum_{(K,L) \in \cal E} \sqrt{T_{K|L}} \big( u_L^{(n+1)} - u_K^{(n+1)} \big)\cdot \sqrt{T_{K|L}}\big( w_L^{(n(t)+1)} - w_K^{(n(t)+1)}\big)\\
      \leqslant \frac{a}{2 \alpha} \sum_{n=n(t)+1}^{n(t+\tau)} \dt^{(n)} \sum_{(K,L) \in \cal E} T_{K|L} \big( u_L^{(n+1)} - u_K^{(n+1)} \big)^2\\
      + \frac{a}{2 \alpha} \sum_{n=n(t)+1}^{n(t+\tau)} \dt^{(n)} \sum_{(K,L) \in \cal E} T_{K|L} \big( w_L^{(n(t)+1)} - w_K^{(n(t)+1)}\big)^2\\
      \leqslant \frac{a}{2 \alpha} \sum_{n=n(t)+1}^{n(t+\tau)} \dt^{(n)} \sum_{(K,L) \in \cal E} T_{K|L} \big( u_L^{(n+1)} - u_K^{(n+1)} \big)^2\\
      + \frac{a}{\alpha^3} \sum_{n=n(t)+1}^{n(t+\tau)} \dt^{(n)} \sum_{(K,L) \in \cal E} T_{K|L} \big( u_L^{(n(t)+1)} - u_K^{(n(t)+1)}\big)^2\\
      + \frac{a}{\alpha \beta^2} \sum_{n=n(t)+1}^{n(t+\tau)} \dt^{(n)} \sum_{(K,L) \in \cal E} T_{K|L} \big( v_L^{(n(t)+1)} - v_K^{(n(t)+1)}\big)^2,
   \end{gather*}
   where the first inequality follows from the relation $\displaystyle \pm s_1 s_2 \leqslant \frac{1}{2} (s_1^2 + s_2^2)$ and the second one follows from the simple inequality $(s_1 + s_2)^2 \leqslant 2(s_1^2 + s_2^2)$. To conclude the proof we integrate above inequalities over $\mathbb R$ with respect to the time variable $t$. Next we apply Proposition \ref{FVMestttt} (for details see  \cite[Lemma 5.5]{eymard_gutnic}). \hskip 5pt $\blacksquare$ \medskip

Proposition \ref{FVMlempom1} together with Proposition \ref{FVMlempom2} immediately give the following corollary.
\begin{coro}\label{FVMestirc}
   Let the assumptions \textit{1.\ 2.} and \textit{3.}\ of Proposition \ref{FVMlempom1} be satisfied. We set $\displaystyle w_{\cal D} = \frac 1 \alpha u_{\cal D}+ \frac 1 \beta v_{\cal D}$. Then, there exists a constant $C >0$, which is independent of the discretization parameters ${\cal D}$ and of $k$, such that
   \begin{equation}\label{FVMeqsptw}
      \int_{\Omega_\xi\times(0,T)} \big(w_{\cal D}(x+\xi,t) - w_{\cal D}(x,t)\big)^2 \dx \dy{t} \leqslant C |\xi| \big(2 \size({\cal D} \big) +|\xi|),
   \end{equation}
   for all $\xi\in\mathbb R^d$ and $\Omega_\xi=\{x \in \mathbb R^d, [x,x+\xi]
   \subset \Omega\}$. Moreover
   \begin{equation}\label{FVMeqttw}
      \int_{\Omega\times(0,T-\tau)} \big(w_{\cal D}(x,t+\tau) - w_{\cal D}(x,t)\big)^2 \dx \dy{t} \leqslant C \big(\size({\cal D}) +\tau \big),
   \end{equation}
   where $\tau\in(0,T)$.{\hskip 5pt $\blacksquare$}
\end{coro} \medskip

\subsection{The limit as $size(\cal D)$ tends to zero and $k$ tends to infinity}

We state below the main convergence results of this paper, first only letting the size of the volume elements and the time steps tend to zero, and then also letting the kinetic rate tend to infinity.
\begin{prop}\label{FVMestcvg}
   We suppose that the hypotheses $\mathcal H$ are satisfied. Let $k>0$ be arbitrary and let $(u^k,\,v^k) \in \Big(L^2\big(0,T;\, H^1 ( \Omega) \big) \Big)^2$ be the unique classical solution of \Pp. Then the sequence $(u^k_{{\cal D}},\, v^k_{{\cal D}})$ of approximate solutions of \Pp\ given by \eqref{FVMschini}, \eqref{FVMscheme} and \eqref{FVMuapp} converges strongly in $\Big(L^2\big(Q_T\big) \Big)^2$ to $(u^k,\, v^k)$ as $\size({\cal D}) \rightarrow 0$. In particular, the sequence $\displaystyle w^k_{{\cal D}} = \frac 1 \alpha u^k_{{\cal D}}+ \frac 1 \beta v^k_{{\cal D}}$ converges strongly in $L^2\big(Q_T\big)$ to the function $\displaystyle w^k = \frac 1 \alpha u^k+ \frac 1 \beta v^k$ as $\size({\cal D}) \rightarrow 0$ and there exist positive constants $C_1,\, C_2$ which do not depend on $k$, such that
   \begin{gather}
      k \int_{\Omega \times(0,T)} \big(r_A(u^k) - r_B(v^k)\big)^2 \dx \leqslant C_1,\label{FVMwarunek}\\
      \int_{\Omega_\xi\times(0,T)} \big(w^k(x+\xi,t) - w^k(x,t)\big)^2 \dx \dy{t} \leqslant
      C_2 |\xi|^2,\label{FVMeqsptwlim}\\
      \int_{\Omega\times(0,T-\tau)} \big(w^k(x,t+\tau) - w^k(x,t)\big)^2 \dx \dy{t} \leqslant
      C_2 \tau \label{FVMeqttw2},
   \end{gather}
   where $\tau\in(0,T),\ \xi\in\mathbb R^d$ and $\Omega_\xi=\{x \in \mathbb R^d, [x,x+\xi]
   \subset \Omega\}$.
\end{prop}
{\hskip -\parindent \bf Proof \hskip 5pt}
   To prove the result we use Corollary \ref{FVMestirc}. The method of proof is similar to that of Theorem \ref{FVMfixkconvth}.
{\hskip 5pt $\blacksquare$}

It is now possible to pass to the limit as $k\rightarrow+\infty$.
\begin{theo}\label{FVMestcvg2}
   Let $\big( u_{\mathcal D}^{k},\,v_{\mathcal D}^{k}\big)$ be the sequence of approximate solutions of \Pp, defined by \eqref{FVMschini}, \eqref{FVMscheme} and \eqref{FVMuapp}. Then
   \begin{equation*}
      u^k_{\mathcal D} \longrightarrow \bigg( \frac{1}{\alpha} \text{id} + \frac{1}{\beta} \eta \bigg)^{-1}\!\!\!(w)
   \end{equation*}
   and
   \begin{equation*}
      v^k_{\mathcal D} \longrightarrow \eta \circ \bigg( \frac{1}{\alpha} \text{id} + \frac{1}{\beta} \eta \bigg)^{-1}\!\!\!(w)
   \end{equation*}
   as  $\size(\mathcal D) \rightarrow 0$ and $k$ tends to $\infty$, where $\eta = r_B^{-1} \circ r_A$ and where $w$ is the unique weak solution of the problem \eqref{FVMNonlDiff1}\,--\,\eqref{FVMNonlDiff2}.
\end{theo}
{\hskip -\parindent \bf Proof \hskip 5pt}
   Let $\displaystyle w^k_{{\cal D}} = \frac{1}{\alpha} u^k_{{\cal D}} + \frac{1}{\beta} v^k_{{\cal D}}$. The estimates from Corollary \ref{FVMestirc}, which are uniform with respect to $k$, permit to apply Proposition \ref{FVMconskol}. As a consequence we deduce the relative compactness in $L^2(Q_T)$ of the sequence $\{ w^k_{{\cal D}} \}$. Then there exist a function $w \in L^2(Q_T)$ and a subsequence $\{w^{k_i}_{{\cal D}_{m}}\}$ such that  $w^{k_i}_{{\cal D}_{m}}$ converges to $w$ strongly in $L^2(Q_T)$ as $k_i$ tends to infinity and $size({\cal D}_{m})$ tends to zero. Theorem \ref{FVMlinf} implies that $w$ is nonnegative and bounded in $Q_T$.
   The inequality \eqref{FVMeqchose}, namely
   \begin{equation*}
      {k_i} \big \| r_A(u^{k_i}_{{\cal D}_m}) - r_B(v^{k_i}_{{\cal D}_m}) \big \|^2_{L^2(Q_T)} \leqslant C
   \end{equation*}
   where the positive constant $C$ is independent of $k_i$ and $size({\cal D}_{m})$, implies that
   \begin{equation*}
      r_A(u^{k_i}_{{\cal D}_m}) - r_B(v^{k_i}_{{\cal D}_m}) \rightarrow 0 \quad \text{in} \quad L^2(Q_T),
   \end{equation*}
   and consequently almost everywhere, as $k_i$ tends to infinity. Then
   \begin{equation*}
      v^{k_i}_{{\cal D}_m} = \eta(u^{k_i}_{{\cal D}_m}) + e^{k_i}_{{\cal D}_m},
   \end{equation*}
   where $\eta(s) = r^{-1}_B \big(r_A(s)\big)$ and $e^{k_i}_{{\cal D}_m}$ tends to zero almost everywhere as $size({\cal D}_m)$ tends to zero and $k_i$ tends to infinity. In view of the hypotheses $\mathcal H\, 4$ the function $\eta(s)$ is well defined on $[0,\infty)$. Moreover,
   \begin{equation*}
      H(u^{k_i}_{{\cal D}_m}) = w^{k_i}_{{\cal D}_m} - \frac{1}{\beta} e^{k_i}_{{\cal D}_m} \rightarrow w \quad \text{a.e.\ in} \quad Q_T,
   \end{equation*}
   where $\displaystyle H(s) = \frac{1}{\alpha} s + \frac{1}{\beta} \eta(s)$. Hypotheses $\mathcal H\, 4$ ensures that the function $H(s)$ has the continuous inverse function. Then Lebesgue's dominated convergence theorem implies that there exists a function pair $(\tilde u,\,\tilde v) \in \big(L^2(Q_T)\big)^2$
   \begin{equation}\label{FVMLebesgue}
      u^{k_i}_{{\cal D}_m} \rightarrow \tilde{u} \quad \text{and} \quad v^{k_i}_{{\cal D}_m} \rightarrow \tilde{v} \quad \text{in} \quad L^2(Q_T)
   \end{equation}
   as $size({\cal D}_m)$ tends to zero and $k_i$ tends to infinity. \\
   Next we identify the limit pair $(\tilde u,\, \tilde v)$. Let $\overline{w}$ be the solution of the problem \eqref{FVMNonlDiff1}\,--\,\eqref{FVMNonlDiff2} and the functions $\overline u$ and $\overline v$ are defined as in \eqref{FVMrepresUV}, namely
   \begin{equation*}
      \overline{u} = \bigg( \frac{1}{\alpha} \text{id} + \frac{1}{\beta} \eta \bigg)^{-1} \hskip -4pt \big( \overline{w} \big),\quad \overline{v} = \eta \circ \bigg( \frac{1}{\alpha} \text{id} + \frac{1}{\beta} \eta \bigg)^{-1} \hskip -4pt \big( \overline{w} \big).
   \end{equation*}
   Let $\eps > 0$ be arbitrary. We have that
   \begin{equation*}
      \rule{20pt}{0pt} \big \| \overline u - \tilde u \big \|_{L^2(Q_T)} \leqslant \big \| \overline u - u^{k_i} \big \|_{L^2(Q_T)} + \big \| u^{k_i} - u_{\mathcal D_m}^{k_i} \big \|_{L^2(Q_T)} + \big \| u_{\mathcal D_m}^{k_i} - \tilde u \big \|_{L^2(Q_T)} \rule{5pt}{0pt}.
   \end{equation*}
   From \eqref{FVMLebesgue} we deduce that there exists $k_0$ and $\delta_0$ such that, for all $k_i \geq k_0$ and all $size({\cal D}_m) \leq \delta_0$, the last term of the inequality above is less then $\eps/3$. From Theorem 1 in \cite{bothe_hilhorst} there exists some ${{\bar k}_0}$, such that for all ${k_i} > {\bar k}_0$, $\big \| \overline u - u^{k_i} \big \|_{L^2(Q_T)} \leqslant \eps/3$. Then, fixing $k_i = \max (k_0,{\bar k}_0)$, we can take $size({\cal D}_m) \leq \delta_0$ small enough so that by  Proposition \ref{FVMestcvg} $\big \| u^{k_i} - u_{\mathcal D_m}^{k_i} \big \|_{L^2(Q_T)} \leqslant \eps/3$.\\
   Since the argument for the $v$-component is similar, this completes the proof.
{\hskip 5pt $\blacksquare$}

\section{Numerical example}\label{FVMnum_resul}

In this section we give an example of an application of the finite volume scheme \eqref{FVMscheme} in one space dimension. \\
For the numerical experiments we choose the reaction of the re\-versible dimerisation of $o$-phenylene\-dio\-xydi\-me\-thyl\-si\-la\-ne ($2,2$-di\-me\-thyl-$1,2,3$-ben\-zo\-dio\-xa\-si\-lo\-le) which has been studied by $^1$\hskip -1pt H NMR spe\-c\-tro\-sco\-py. The kinetics of this reaction can be described quantitatively by a bimolecular lO-ring formation reaction and a mono\-mo\-le\-cular back\-re\-a\-ction (for further details we refer to  Meyer, Klein and Weiss \cite{meyer_klein_weiss}). \\
Since the reaction is of the type $2 \mathcal A \rightleftharpoons \mathcal B$, the reaction terms take the form
\begin{equation*}
   r_A(u) = k_1 u^2 \quad \text{and} \quad r_B(v) = k_2 v.
\end{equation*}
Moreover $\alpha = 2$ and $\beta = 1$. For this particular process benzene was chosen as a solvent. Then it was possible to estimate rate constants for both reactions at the temperature $T = 298 K$,
\begin{equation*}
   k_1 \approx 1,072 \cdot 10^{-4} \text{L}^2 \text{mol}^{-2} \quad \text{and} \quad k_2 \approx 2,363 \cdot 10^{-6} \text{L}^2 \text{mol}^{-2}
\end{equation*}
and diffusion coefficients
\begin{equation*}
   a \approx 1,579 \cdot 10^{-9} \text{m}^2 \text{s}^{-1} \quad \text{and} \quad b \approx 1,042 \cdot 10^{-9} \text{m}^2 \text{s}^{-1}.
\end{equation*}
In the experiment we set $k=1$ for the chemical kinetics factor. We remark that it is equivalent to the situation when coefficients $a,\ b,\ k_1$ and $k_2$ are of order $1$ and $k$ is of order $10^{4}$. In fact, we can multiply the system \eqref{FVMmain_eq} by $10^{9}$ and change the time scale as $t \mapsto 10^9 t$. However the above reasoning is formally correct and shows in an explicit way the order of the kinetics factor $k$; in our example we decided to keep constants in the form given by the spectroscopic analysis.

Figure 1 shows the initial conditions $u_0(x)$ and $v_0(x)$, defined as
\begin{equation}\label{FVM_init_num_def_u}
   u_0(x) = \begin{cases}
      0 & \text{ for } x \in [0,0.03]\\
      \frac{1}{2} \sin \Big( \frac{50 \pi}{7}(x - 0.03) \Big) & \text{ for } x \in [0.03,0.1],
   \end{cases}
\end{equation}
and
\begin{equation}\label{FVM_init_num_def_v}
   v_0(x) = \begin{cases}
      \frac{1}{4} \cos \Big( \frac{50 \pi}{7}x \Big) & \text{ for } x \in [0,0.07],\\
      0 & \text{ for } x \in [0.07,0.1]
   \end{cases}
\end{equation}
On Figures 3 and 5 we see the time evolution of the solution $\big( u_{{\cal D}}^k,v_{{\cal D}}^k \big)$ until the times $T_\text{max}^1 = 10^{5}\ [s]$ and $T_\text{max}^2 = 10^{11}\ [s]$, respectively. Then we follow the evolution of the solution $w_{{\cal D}}$ of the nonlinear diffusion problem \eqref{FVMNonlDiff1}\,--\,\eqref{FVMNonlDiff2} for initial condition deduced from that used in the reaction\,--\,diffusion \Pp. We have used a uniform mesh with $h = 0,002$ and initial time step $t_\delta = 10^{-8}$ to obtain the approximate solution $\big( u_{{\cal D}}^k,v_{{\cal D}}^k \big)$ and $t_\delta = 10^{-6}$ to obtain the approximate solution $w_{{\cal D}}$to the nonlinear diffusion problem.

We can use the approximate solution $w_{{\cal D}}$ and the formulas \eqref{FVMrepresUV} to define functions $u_{{\cal D}}^w$ and $v_{{\cal D}}^w$. Indeed, let
\begin{equation*}
   u_{{\cal D}}^w = h \big( w_{{\cal D}} \big) \quad \text{and} \quad v_{{\cal D}}^w = g \big( h \big( w_{{\cal D}} \big) \big),
\end{equation*}
where
\begin{equation}\label{FVMdef_h}
   h(y) = \frac{1}{2} \bigg( \sqrt{\Big( \frac{\alpha k_1}{\beta k_2} \Big)^2 + y\,\frac{4 k_2}{\beta k_1}} - \frac{\alpha k_2}{\beta k_1}\bigg)
\end{equation}
and
\begin{equation}\label{FVMdef_g}
   g(h) = h\,\frac{a}{\alpha} + h^2\,\frac{b k_1}{\beta k_2} \hskip 3pt.
\end{equation}

Proceeding in the similar way as in the proof of Theorem \ref{FVMestcvg2}, we write for the $u$-component that
\begin{multline}\label{FVMtrobo}
   \big \| u_{{\cal D}}^k - h(w_{{\cal D}})\big \|_{L^2(Q_T)} \leqslant \big \| u_{{\cal D}}^k - u^k \big \|_{L^2(Q_T)} \\
   + \big \| u^k - h(w) \big \|_{L^2(Q_T)} + \big \| h(w) - h(w_{{\cal D}}) \big \|_{L^2(Q_T)}.
\end{multline}
We simultaneously pass to the limit as $\size({\cal D}) \rightarrow 0$ and $k\rightarrow \infty$. From Theorem \ref{FVMestcvg2} we immediately deduce that the first term on the right hand side of \eqref{FVMtrobo} tends to zero. The same conclusion holds for the two other terms. Indeed, \cite[Theorem 1, Sec.\ 3]{bothe_hilhorst} implies that $\big \| u^k - h(w) \big \|_{L^2(Q_T)}$ tends to zero as $k \rightarrow \infty$. Moreover, \cite[Theorem 5.1]{eymard_slimane} yields that $\big \| w - w_{{\cal D}} \big \|_{L^2(Q_T)} \rightarrow 0$ as $ \size({\cal D}) \rightarrow 0$ and since the function $h$ is well defined and continuous, we conclude that for every small $\eps >0$ there exist ${\cal D}$ small enough and $k$ large enough so that
\begin{equation*}
   \big \| u_{{\cal D}}^k - h(w_{{\cal D}})\big \|_{L^2(Q_T)} \leqslant \eps.
\end{equation*}
We proceed in the same way to show that for every small $\eps >0$ there exist ${\cal D}$ small enough and $k$ large enough so that
\begin{equation*}
   \big \| v_{{\cal D}}^k - g(h(w_{{\cal D}}))\big \|_{L^2(Q_T)} \leqslant \eps.
\end{equation*}

The results from our numerical experiment agree with above analysis, since
\begin{align*}
   & \max_{x \in \Omega} \Big| u_{{\cal D}}^k(x,T_\text{max}^1) - h(w_{{\cal D}})(x, T_\text{max}^1) \Big|_{\infty} \simeq 4,74 \cdot 10^{-3}, \\
   & \max_{x \in \Omega} \Big| v_{{\cal D}}^k(x,T_\text{max}^1) - g(h(w_{{\cal D}})(x, T_\text{max}^1)) \Big|_{\infty} \simeq 4,032 \cdot 10^{-3},
\end{align*}
whereas
\begin{align*}
   & \max_{x \in \Omega} \Big| u_{{\cal D}}^k(x,T_\text{max}^2) - h(w_{{\cal D}})(x, T_\text{max}^2) \Big|_{\infty} \simeq 3,121 \cdot 10^{-14}, \\
   & \max_{x \in \Omega} \Big| v_{{\cal D}}^k(x,T_\text{max}^2) - g(h(w_{{\cal D}})(x, T_\text{max}^2)) \Big|_{\infty} \simeq 1,84 \cdot 10^{-13}.
\end{align*}


\hskip -\parindent \begin{minipage}{0.45\textwidth}
   \begin{center}
      \includegraphics[width = \textwidth]{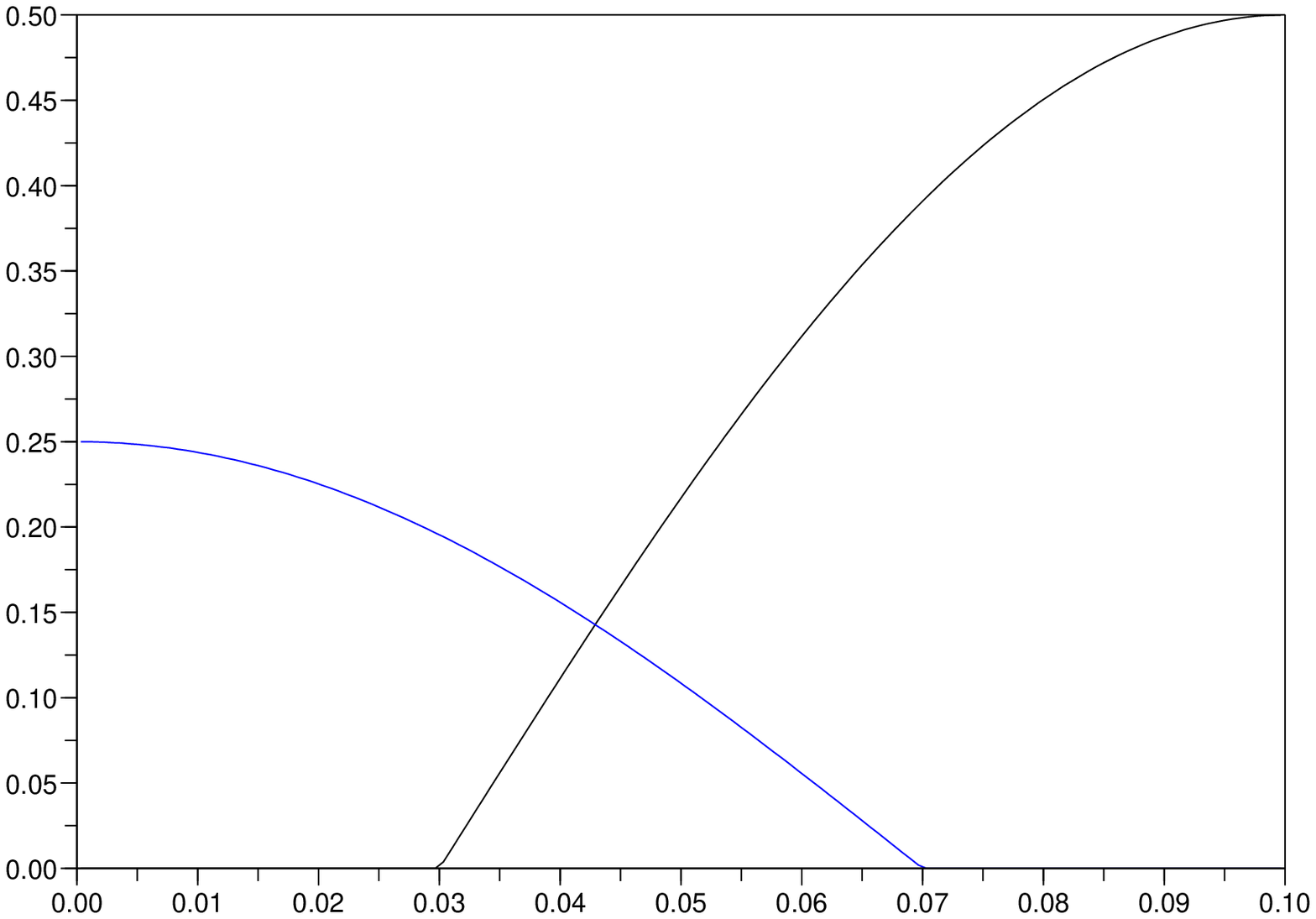}\\
      \textbf{Figure 1:} Initial conditions $u_0$ and $v_0$, defined in \eqref{FVM_init_num_def_u} and \eqref{FVM_init_num_def_v}.
   \end{center}
\end{minipage}
\begin{minipage}{0.45\textwidth}
   \begin{center}
      \includegraphics[width = \textwidth]{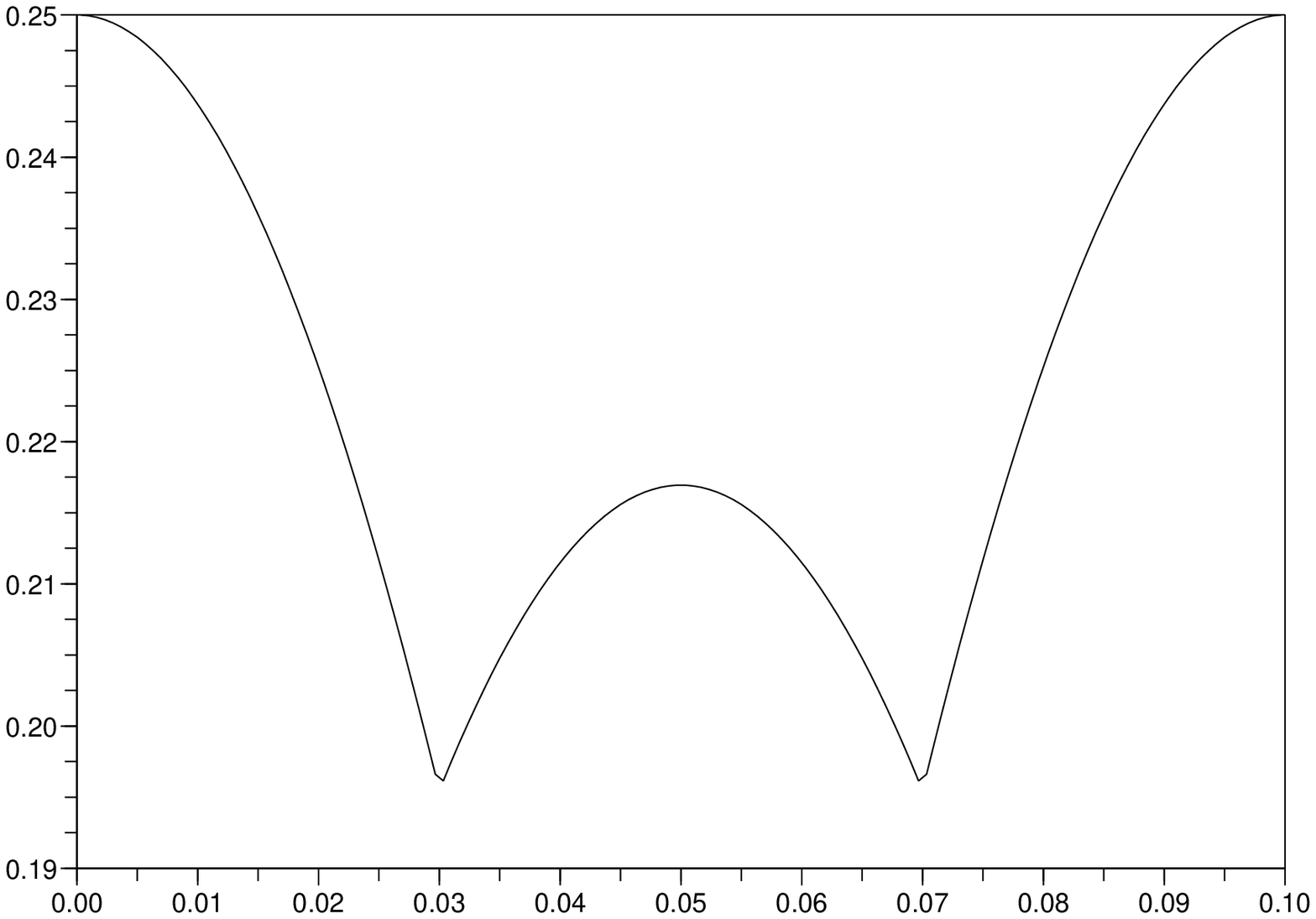}\\
      \textbf{Figure 2:} Initial condition $w_0 = \frac{1}{\alpha}u_0 + \frac{1}{\beta} v_0$, where $u_0,\ v_0$ are defined in \eqref{FVM_init_num_def_u} and \eqref{FVM_init_num_def_v}.
   \end{center}
\end{minipage}
\begin{minipage}{0.45\textwidth}
   \begin{center}
      \includegraphics[width = \textwidth]{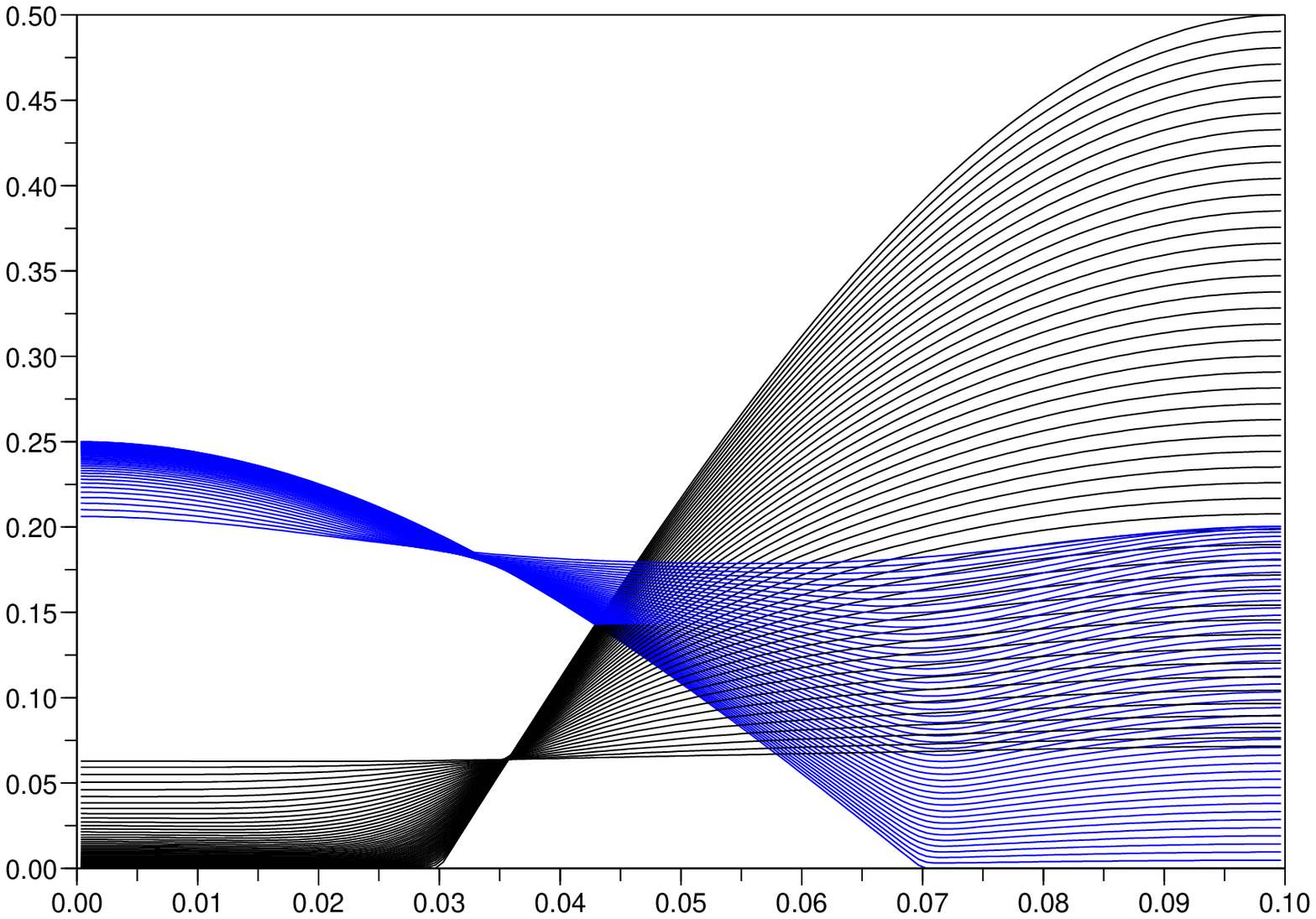}\\
      \textbf{Figure 3:} Evolution of $u_{\text{ap}}$ and $v_{\text{ap}}$ on the time interval $(0,\ldots,T_\text{max}^1)$.
   \end{center}
\end{minipage}
\begin{minipage}{0.45\textwidth}
   \begin{center}
      \includegraphics[width = \textwidth]{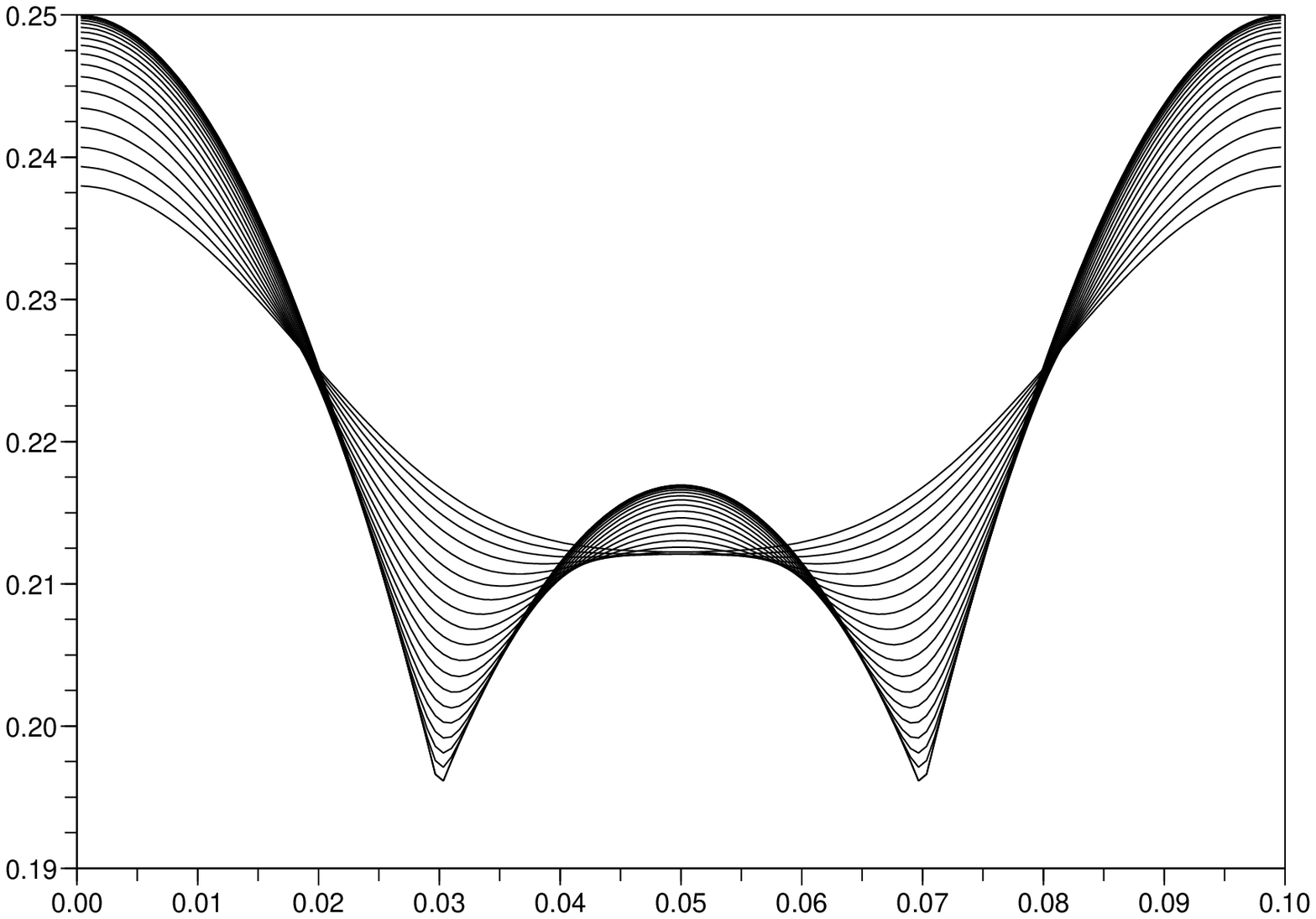}\\
      \textbf{Figure 4:} Evolution of $w_{\text{ap}}$ on the time interval $(0,\ldots,T_\text{max}^1)$.
   \end{center}
\end{minipage}
\begin{minipage}{0.45\textwidth}
   \begin{center}
      \includegraphics[width = \textwidth]{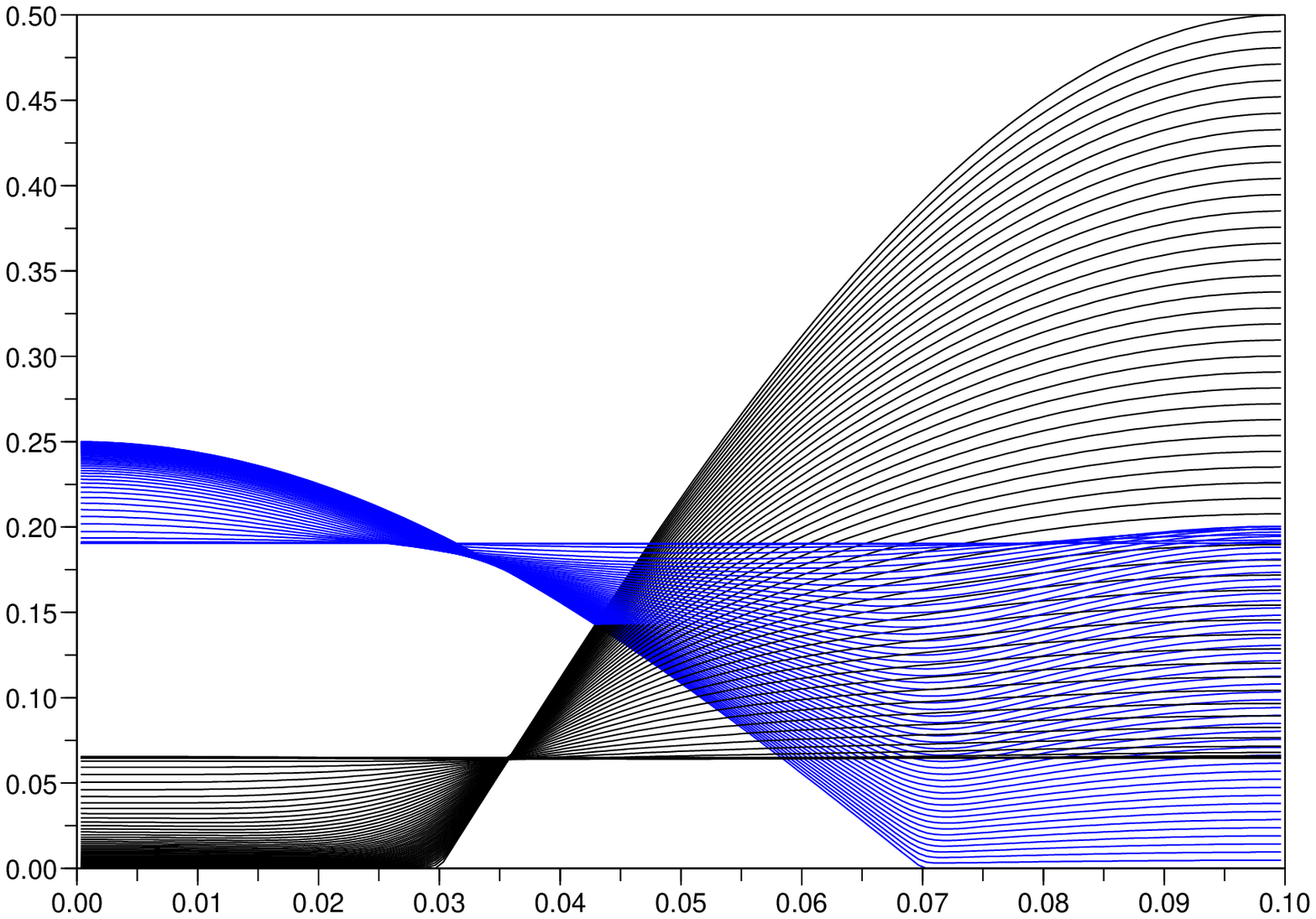}\\
      \textbf{Figure 5:} Evolution of $u_{\text{ap}}$ and $v_{\text{ap}}$ on the time interval $(0,\ldots,T_\text{max}^2)$.
   \end{center}
\end{minipage}
\begin{minipage}{0.45\textwidth}
   \begin{center}
      \includegraphics[width = \textwidth]{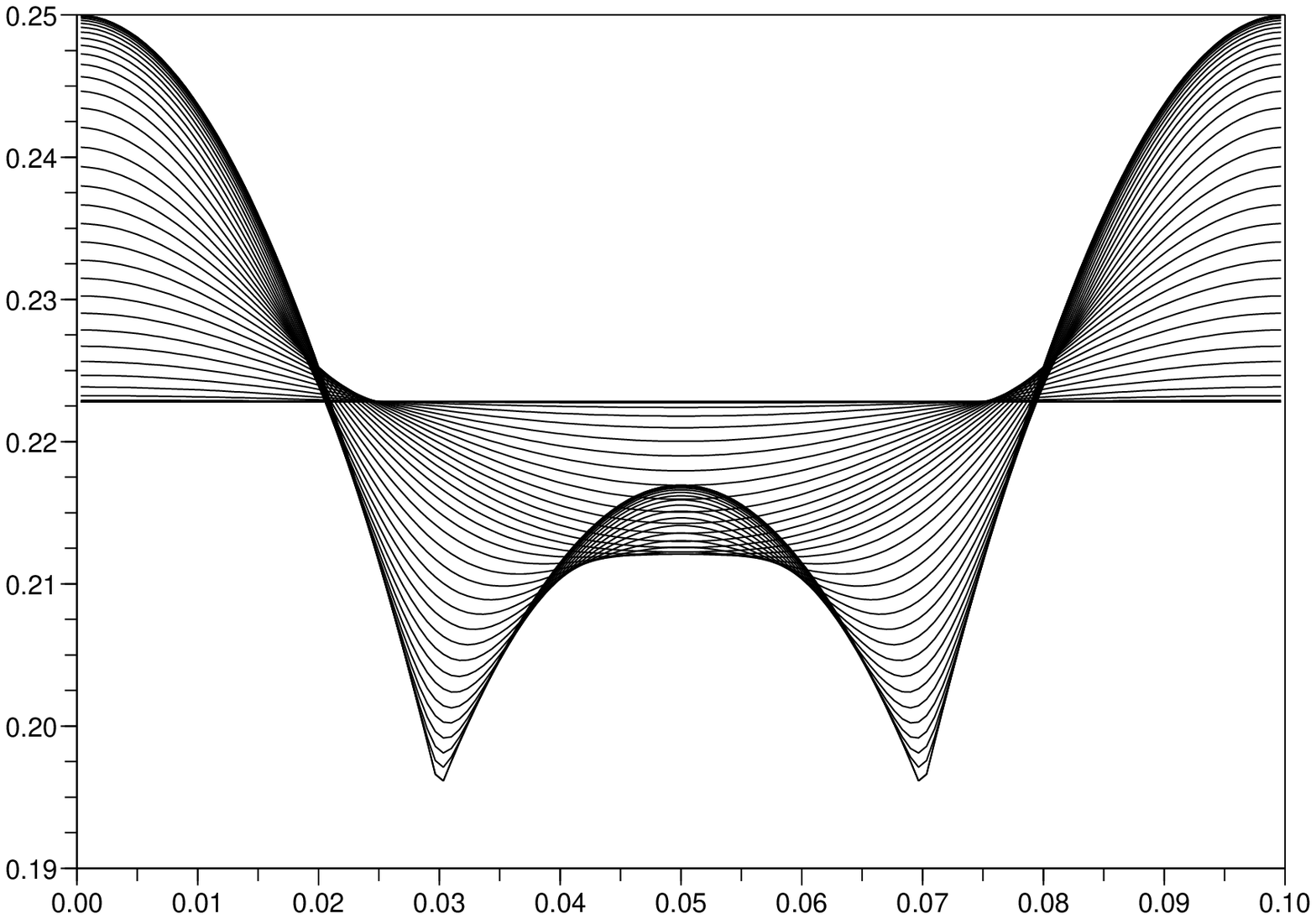}\\
      \textbf{Figure 6:} Evolution of $w_{\text{ap}}$ on the time interval $(0,\ldots,T_\text{max}^2)$.
   \end{center}
\end{minipage}

\vskip 30pt

In order to show the accuracy of our method in the fast reaction limit, let us increase the kinetics parameter $k$ in \Pp\, keeping all other data as previously. Let
\begin{equation*}
   \mathcal J_u^k = \max_{x \in \Omega} \Big| u_{{\cal D}}^k(x,T_\text{max}^2) - h(w_{{\cal D}})(x, T_\text{max}^2) \Big|_{\infty}
\end{equation*}
and
\begin{equation*}
   \mathcal J_v^k = \max_{x \in \Omega} \Big| v_{{\cal D}}^k(x,T_\text{max}^2) - g(h(w_{{\cal D}})(x, T_\text{max}^2)) \Big|_{\infty}.
\end{equation*}

Table 1 below shows the numerical results.

\bigskip

\begin{center}
   \begin{tabular}{c r@{$\cdot$}l r@{.}l}
      $k$ & \multicolumn{2}{c}{$\mathcal J_u^k$} & \multicolumn{2}{c}{$\mathcal J_v^k$} \\
      \hline \hline
      \rule{0pt}{12pt}
      $10^{-7}$ & $2,3498$&$10^{-2}$ &  $4,5272$&$10^{-2}$ \\
      $10^{-6}$ & $3,976$&$10^{-3}$ &   $1,988$&$10^{-3}$ \\
      $10^{-5}$ & $4,699$&$10^{-8}$ &   $2,3498$&$10^{-8}$ \\
      $10^{-4}$ & $9,3725$&$10^{-10}$ & $4,6862$&$10^{-10}$ \\
      $10^{-3}$ & $2,498$&$10^{-16}$ &  $9,992$&$10^{-16}$ \\
      $10^{-2}$ & $1,7892$&$10^{-12}$ & $8,95$&$10^{-13}$ \\
      $10^{-1}$ & $4,163$&$10^{-17}$ &  $9,992$&$10^{-16}$ \\
      $10^{0}$  & $3,121$&$10^{-17}$ &  $1,84$&$10^{-16}$ \\
   \end{tabular}\\
   \rule{0pt}{20pt}\textbf{Table 1:} The accuracy of our method in the fast reaction limit, when the kinetics parameter $k$ in \Pp\,  increases and all other data are unchanged.
\end{center}

\section{Appendix}

The proof of the following result can be found in \cite[Lemma 5.3 and Lemma 5.4]{eymard_gutnic}.
\begin{lemm}\label{FVMestttt}
   Let $(t^{(n)})_{n\in{\mathbb Z}}$ be a strictly increasing sequence of real numbers
   such that $\lim_{n\rightarrow-\infty} t^{(n)} = -\infty$ and $\lim_{n\rightarrow\infty} t^{(n)} = \infty$. Moreover, let $\dtn := t^{(n+1)} - t^{(n)}$ be uniformly bounded. For all $t\in\mathbb R$ we denote by $n(t)$ an integer $n$, such that $t\in[t^{(n)},t^{(n+1)})$. Let $(a^{(n)})_{n\in{\mathbb Z}}$ be a family of nonnegative real values such that $a^{(n)} \neq 0$ for finitely many $n \in \mathbb Z$. Then, for all $\tau\in (0,+\infty)$ and $\zeta\in\mathbb R$
   \begin{gather}
      \label{FVMesttt1} \int_{\mathbb R} \sum_{n=n(t)+1}^{n(t+\tau)} \big(\dtn a^{(n+1)}\big) \dy{t} = \tau \sum_{n\in{\mathbb Z}} \big(\dtn a^{(n+1)}\big), \\
      \label{FVMesttt2} \int_{\mathbb R} \bigg(\sum_{n=n(t)+1}^{n(t+\tau)}\dtn \bigg) a^{(n(t+\zeta)+1)} \dy{t} \leqslant \big(\tau + \max_{n\in{\mathbb Z}} \dtn\big) \sum_{n\in{\mathbb Z}} \big(\dtn a^{(n+1)}\big).
   \end{gather}
\end{lemm}

\medskip

The following proposition is a direct corollary from Fr\'echet--Kolmogorov Theorem (see Brezis \cite[Theorem IV.25, page 72]{brezis}).
\begin{prop}\label{FVMconskol}
   Let $\mathcal O$ be a bounded and open subset of $\mathbb R^{d+1}$, $d = 1,2$ or $3$. Let $(w_n)_{n\in\N}$  be a sequence of functions $w_n(x,t):\ \mathbb R^d \times \mathbb R \rightarrow \mathbb R$, such that
   \begin{enumerate}
      \item for all $n\in\N$, $w_n\in L^\infty(\mathcal O)$ and there exists a constant $C_b >0$ which does not depend on $n$, such that $\| w_n \|_{L^\infty(\mathcal O)} \leqslant C_b$,
      \item there exist positive constants $C_1,\, C_2$ and a sequence of nonnegative real values $(\mu_n)_{n\in\N}$, such that $\lim_{n\rightarrow\infty} \mu_n = 0$ and
      \begin{equation*}
         \int_{\mathcal O_{\xi,\tau}} \big(w_n(x + \xi,t + \tau) - w_n(x,t)\big)^2 \dx \dy{t} \leqslant C_1 |\xi| + C_2 \tau + \mu_n,
      \end{equation*}
      for $\xi\in\mathbb R^d,\ \tau \in \mathbb R,\ n\in\on$ and
      \begin{equation*}
         \rule{20pt}{0pt}\mathcal O_{\xi,\tau} = \big \{ (x,t)\in \mathbb R^{d+1}:\text{\ the interval\ } \big[(x,t),(x+\xi,t+\tau)\big] \text{\ lies in\ } {\mathcal O} \big \}.\rule{20pt}{0pt}
      \end{equation*}
   \end{enumerate}
   Then there exists a subsequence of $(w_n)_{n\in\N}$, denoted again by $(w_n)_{n\in\N}$ and a function $w \in L^\infty(\mathcal O)$ such that $w_n \rightarrow w$ in $L^2(\mathcal O)$, as $n\rightarrow\infty$.
\end{prop}

\medskip

\begin{lemm}\label{FVMoszacLN}
   Let $A>0$ and a function $r(s)$ satisfying
   \begin{equation}\label{FVMhypsupLEM}
   \begin{split}
      & r \in C^1(\mathbb R),\ r'(\cdot) > 0 \text{\quad on \quad}(0,+\infty),\\
      & r(0) = 0, \text{\quad and \quad} \limsup_{s \rightarrow 0^+} \frac{s r'(s)}{r(s)} < \infty
   \end{split}
   \end{equation}
   be given.\\
   Then there exists $\eps_0 > 0$ only depending on $r(s)$ and $C > 0$, depending only on $r(s)$ and $A$, such that for all $\eps \in (0,\eps_0)$ and for all $u \in [0,A]$ the inequality
   \begin{equation*}
      \big| \ln r(u + \eps) \big| \leqslant C \big( |\ln \eps| + 1 \big)
   \end{equation*}
   holds. 
\end{lemm}
{\hskip -\parindent \bf Proof \hskip 5pt}
   Let $\alpha = \displaystyle \limsup_{s \rightarrow 0^+} \frac{s r'(s)}{r(s)} +1$. There exists a constant $\eps_0 > 0$, such that for all $s \in (0,\eps_0)$ 
   \begin{equation}\label{FVMineqalpha}
      \frac{r'(s)}{r(s)} \leqslant \frac{\alpha}{s}.
   \end{equation}
   Let $\eps \in (0,\eps_0)$ and $u \in [0,A]$. Then, either $r(u+\eps) \geqslant 1$, which implies
   \begin{equation*}
      \ln r(u+\eps) \leqslant \ln r(A+\eps_0)
   \end{equation*}
   or $r(u+\eps) \leqslant 1$. In that case
   \begin{equation*}
      \big| \ln r(u+\eps) \big| \leqslant \big| \ln r(\eps) \big|,
   \end{equation*}
   holds. Integrating inequality \eqref{FVMineqalpha} over the interval $[\eps,\eps_0]$, we obtain
   \begin{equation*}
      \big| \ln r(\eps) \big| \leqslant \big| \ln r(\eps_0) \big| + \alpha \big( |\ln \eps_0| + |\ln \eps| \big).
   \end{equation*}
   Setting $C = \max \big \{ \alpha,\, \alpha |\ln \eps_0| + |\ln r(\eps_0)|,\, \ln \big (r(A+\eps_0)\big) \big \}$ we conclude the proof.
{\hskip 5pt $\blacksquare$}


\end{document}